\documentclass[onefignum,onetabnum]{siamart190516}

\usepackage{lipsum}
\usepackage{amsfonts}
\usepackage{graphicx}
\usepackage{epstopdf}
\usepackage{algorithmic}
\ifpdf
  \DeclareGraphicsExtensions{.eps,.pdf,.png,.jpg}
\else
  \DeclareGraphicsExtensions{.eps}
\fi

\usepackage{caption}
\usepackage{cite}
\usepackage{float}
\floatstyle{plaintop}
\restylefloat{table}

\def\Ker{{\rm Ker \,}}

\def\Def{\stackrel{\mathrm{def}}{=}}

\def\beq{\begin{equation}}
\def\eeq{\end{equation}}

\def\R{\mathbb{R}}
\def\E{\mathbb{E}}

\def\BI{\begin{itemize}}
\def\EI{\end{itemize}}
\def\II{\item}

\newcommand{\refLE}[1]{\ensuremath{\stackrel{(\ref{#1})}{\leq}}}
\newcommand{\refEQ}[1]{\ensuremath{\stackrel{(\ref{#1})}{=}}}
\newcommand{\refGE}[1]{\ensuremath{\stackrel{(\ref{#1})}{\geq}}}

\newcommand{\refPE}[1]{\ensuremath{\stackrel{(\ref{#1})}{\preceq}}}
\newcommand{\refSE}[1]{\ensuremath{\stackrel{(\ref{#1})}{\succeq}}}
\newcommand{\refIN}[1]{\ensuremath{\stackrel{(\ref{#1})}{\in}}}

\DeclareMathOperator*{\argmax}{argmax}

\def\Tr{{\rm Trace}}
\def\DFP{{\rm DFP}}
\def\BFGS{{\rm BFGS}}
\def\SR{{\rm SR1}}
\def\Broyd{{\rm Broyd}}
\def\Unif{{\rm Unif}}

\newsiamremark{remark}{Remark}
\newsiamremark{hypothesis}{Hypothesis}
\newsiamremark{example}{Example}
\crefname{hypothesis}{Hypothesis}{Hypotheses}
\newsiamthm{claim}{Claim}

\renewcommand\arraystretch{1.5}
\def\ba{\begin{array}}
\def\ea{\end{array}}
\def\beann{\begin{eqnarray*}}
\def\eeann{\end{eqnarray*}}
\def\bea{\begin{eqnarray}}
\def\eea{\end{eqnarray}}
\def\bal{\renewcommand\arraystretch{2}\begin{tabular}{|l|}}
\def\eal{\end{tabular}}

\def\BT{\begin{theorem}}
\def\ET{\end{theorem}}
\def\BL{\begin{lemma}}
\def\EL{\end{lemma}}
\def\BC{\begin{corollary}}
\def\EC{\end{corollary}}
\def\BE{\begin{example}}
\def\EE{\end{example}}
\def\BD{\begin{definition}}
\def\ED{\end{definition}}
\def\BR{\begin{remark}}
\def\ER{\end{remark}}
\def\BAS{\begin{assumption}}
\def\EAS{\end{assumption}}
\def\BI{\begin{itemize}}
\def\EI{\end{itemize}}
\def\BCA{\begin{cases}}
\def\ECA{\end{cases}}

\def\BMP{\begin{minipage}{9.5cm}}
\def\EMP{\end{minipage}}
\def\MPT{\begin{minipage}{11.5cm}}
\def\EPT{\end{minipage}}

\def\la{\langle}
\def\ra{\rangle}

\headers{Greedy Quasi-Newton Methods}{%
   Anton Rodomanov and Yurii Nesterov}

\title{Greedy Quasi-Newton Methods\\
with Explicit Superlinear Convergence\thanks{Received by the editors February
20, 2020; accepted for publication (in revised form) November 18, 2020;
published electronically March 1, 2021.}
\URL{10.1137/20M1320651}
\funding{The research results of this paper were obtained
with support of ERC Advanced Grant 788368.}}

\author{Anton Rodomanov\thanks{Institute of Information and
Communication Technologies, Electronics and Applied
Mathematics (ICTEAM), Catholic University of Louvain (UCL), 1348
Louvain-la-Neuve, Belgium (\email{anton.rodomanov@uclouvain.be}).}
\and Yurii Nesterov\thanks{Center for Operations Research
and Econometrics (CORE), Catholic University of Louvain
(UCL), 1348 Louvain-la-Neuve, Belgium (\email{yurii.nesterov@uclouvain.be}).}}

\usepackage{amsopn}

\ifpdf
\hypersetup{
  pdftitle={Greedy Quasi-Newton Methods},
  pdfauthor={A. Rodomanov, Y. Nesterov}
}
\fi

\slugger{siopt}{2021}{31}{1}{785--811}

\begin{document}

\maketitle

\begin{abstract}
In this paper, we study greedy variants of quasi-Newton
methods. They are based on the updating formulas from a
certain subclass of the Broyden family. In particular, this
subclass includes the well-known DFP, BFGS and SR1 updates.
However, in contrast to the classical quasi-Newton methods,
which use the difference of successive iterates for updating
the Hessian approximations, our methods apply basis vectors,
greedily selected so as to maximize a certain measure of
progress. For greedy quasi-Newton methods, we establish an
explicit non-asymptotic bound on their rate of local
superlinear convergence, as applied to minimizing strongly
convex and strongly self-concordant functions (and, in
particular, to strongly convex functions with Lipschitz
continuous Hessian). The established superlinear convergence
rate contains a contraction factor, which depends on the
square of the iteration counter. We also show that greedy
quasi-Newton methods produce Hessian approximations whose
deviation from the exact Hessians linearly converges to
zero.
\end{abstract}



\begin{DOI}
   10.1137/20M1320651
\end{DOI}

\section{Introduction}

\subsection{Motivation}
Quasi-Newton methods have a reputation of the most efficient
numerical schemes for smooth unconstrained optimization. The
main idea of these algorithms is to approximate the standard
Newton method by replacing the exact Hessian with some
approximation, which is updated between iterations according
to special formulas. There exist numerous variants of
quasi-Newton algorithms that differ mainly in the rules of
updating Hessian approximations. The three most popular are
the \emph{Davidon--Fletcher--Powell (DFP)} method
\cite{Davidon1959,FletcherPowell1963}, the
\emph{Broyden--Fletcher--Goldfarb--Shanno (BFGS)} method
\cite{Broyden1970p1,Broyden1970p2,Fletcher1970,Goldfarb1970,Shanno1970},
and the \emph{Symmetric Rank~1 (SR1)} method
\cite{Davidon1959,Broyden1967}. For a general overview of
the topic, see \cite{DennisMore1977} and \cite[Ch.
6]{NocedalWright2006}; also see \cite{LewisOverton2013} for
the application of quasi-Newton methods for non-smooth
optimization.

The most attractive feature of quasi-Newton methods is their
\emph{superlinear convergence}, which was first established
in the 1970s
\cite{Powell1971,BroydenDennisMore1973,DennisMore1974}.
Namely, for several standard quasi-Newton methods (such as
DFP and BFGS), it was proved that the ratio of successive
residuals tends to zero as the number of iterations goes to
infinity. However, the authors did not obtain any explicit
bounds on the corresponding \emph{rate} of this superlinear
convergence. For example, it is unknown whether the
residuals convergence like $O(c^{k^2})$, where $c \in (0,
1)$ is some constant and $k$ is the iteration counter, or
$O(c^{k^3})$, or $O(k^{-k})$, or somehow else. Thus, despite
the qualitative usefulness of the mentioned result, it
still lacks quantitative estimates of the rate of
convergence. Although many other works on quasi-Newton
methods have appeared since then, to our knowledge, all of
them still contain only asymptotic results (see e.g.
\cite{Stachurski1981,GriewankToint1982,ByrdNocedalYuan1987,ByrdNocedal1989,EngelsMartinez1991,YabeYamaki1996,WeiYuYuanLian2004,YabeOgasawaraYoshino2007,MokhtariEisenRibeiro2018,SunTranDinh2019,GaoGoldfarb2019}).
Thus, up to now, there are still no \emph{explicit} and
\emph{non-asymptotic} estimates of the rate of superlinear
convergence of quasi-Newton methods.

In this work, we make a first step towards obtaining such
estimates. We propose new quasi-Newton methods, which are
based on the updating formulas from a certain subclass of
the Broyden family \cite{Broyden1967}. In particular, this
subclass contains the DFP, BFGS and SR1 updates. However, in
contrast to the classical quasi-Newton methods, which use
the difference of successive iterates for updating the
Hessian approximations, our methods apply \emph{basis
vectors}, greedily selected to maximize a certain measure of
progress. For greedy quasi-Newton methods, we establish an
explicit non-asymptotic bound on their rate of local
superlinear convergence, which contains a contraction
factor, depending on the \emph{square} of the iteration
counter. We also show that these methods produce Hessian
approximations whose deviation from the exact Hessians
converges to zero at a \emph{linear} rate. In contrast, it
is known that the standard quasi-Newton methods, in general,
cannot ensure the convergence of the Hessian approximations
to the true Hessian (see e.g.
\cite{DennisMore1974})\footnote{However, it is worth
mentioning that there are some settings, in which the
standard SR1 method indeed yields convergence to the true
Hessian (for more details, see \cite{ConnGouldToint1991}).}.

The idea of using basis vectors in quasi-Newton methods goes
back at least to so-called \emph{methods of dual directions}
\cite{PshenichnyiDanilin1978}, for which it is also possible
to prove both local superlinear convergence of the iterates
and convergence of the Hessian approximations. However,
similarly to the standard quasi-Newton methods, all
corresponding results are only asymptotic. In any case,
despite to the fact that the greedy quasi-Newton methods,
presented in this paper, are based on the same idea, their
construction and analysis are significantly different.

Finally, let us mention that recently there have been
proposed some randomized variants of quasi-Newton
algorithms, which also use nonstandard directions for
updating Hessian approximations
\cite{GowerGoldfarbRichtarik2016,GowerRichtarik2017,KovalevEtAl2020}.

\subsection{Contents}

In Section~\ref{sec-gr-qn}, we discuss a class of
quasi-Newton updating rules for approximating a self-adjoint
positive definite linear operator. We present a special
greedy strategy for selecting an update direction, which
ensures a linear convergence rate in approximating the
target operator. In Section~\ref{sec-quad}, we analyze
greedy quasi-Newton methods, applied to the problem of
minimizing a quadratic function. We show that these methods
have a global linear convergence rate, comparable to that of
the standard gradient descent, and also a superlinear
convergence rate, which contains a contraction factor,
depending on the square of the iteration counter. In
Section~\ref{sec-gen}, we show that similar results also
hold in a more general setting of minimizing a strongly
convex and strongly self-concordant function (and, in
particular, a strongly convex function with Lipschitz
continuous Hessian), provided that the starting point is
chosen sufficiently close to the solution. The main
difficulty here, compared to the quadratic case, is that the
Hessian of the objective function is no longer constant,
resulting in the necessity to apply a special
\emph{correction strategy} to keep the Hessian
approximations under control. Finally, in
Section~\ref{sec-experim}, we present some preliminary
computational results.

\subsection{Notation}

In what follows, $\E$ denotes an arbitrary $n$-dimensional
real vector space. Its dual space, composed of all linear
functionals on $\E$, is denoted by $\E^*$. The value of a
linear function $s \in \E^*$, evaluated at a point $x \in
\E$, is denoted by $\la s, x
\ra$.

For a smooth function $f : \E \to \R$, we denote by $\nabla
f(x)$ and $\nabla^2 f(x)$ its gradient and Hessian
respectively, evaluated at a point $x \in \E$. Note that
$\nabla f(x) \in \E^*$, and $\nabla^2 f(x)$ is a
self-adjoint linear operator from $\E$ to $\E^*$.

The partial ordering of self-adjoint linear operators is
defined in the standard way. We write $A \preceq A_1$
for $A, A_1 : \E \to \E^*$ if $\la (A_1 - A) x, x \ra
\geq 0$ for all $x \in \E$, and $W \preceq W_1$ for
$W, W_1 : \E^* \to \E$ if $\la s, (W_1 - W) s \ra \geq
0$ for all $s \in \E^*$.

Any self-adjoint positive definite linear operator $A : \E
\to \E^*$ induces in the spaces $\E$ and $\E^*$ the
following pair of conjugate Euclidean norms:
\beq\label{def-norms}
\ba{rcll}
\| h \|_A &\Def& \la A h, h \ra^{1/2}, \qquad& h \in \E, \\
\| s \|_A^* &\Def& \la s, A^{-1} s \ra^{1/2}, \qquad& s \in
\E^*.
\ea
\eeq
When $A = \nabla^2 f(x)$, where $f : \E \to \R$ is a smooth
function with positive definite Hessian, and $x \in \E$, we
prefer to use notation $\| \cdot \|_x$ and $\| \cdot
\|_x^*$, provided that there is no ambiguity with the
reference function $f$.

Sometimes, in the formulas, involving products of linear
operators, it is convenient to treat $x \in \E$ as a
linear operator from $\R$ to $\E$, defined by $x \alpha =
\alpha x$, and $x^*$ as a linear operator from $\E^*$ to
$\R$, defined by $x^* s = \la s, x \ra$. In this case, $x
x^*$ is a rank-one self-adjoint linear operator from $\E^*$
to $\E$, acting as follows:
$$
\ba{rcl}
(x x^*) s &=& \la s, x \ra x, \qquad s \in \E^*.
\ea
$$
Likewise, any $s \in \E^*$ can be treated as a linear
operator from $\R$ to $\E^*$, defined by $s \alpha = \alpha
s$, and $s^*$ as a linear operator from $\E$ to $\R$,
defined by $s^* x = \la s, x \ra$. Then, $s s^*$ is a
rank-one self-adjoint linear operator from $\E$ to $\E^*$.

For two self-adjoint linear operators $A : \E \to \E^*$
and $W : \E^* \to \E$, define
$$
\ba{rcl}
\la W, A \ra &\Def& \Tr(W A).
\ea
$$
Note that $W A$ is a linear operator from $\E$ to itself,
and hence its trace is well-defined (it coincides with the
trace of the matrix representation of $W A$ with respect to
an arbitrary chosen basis in the space $\E$, and the result
is independent of the particular choice of the basis).
Observe that $\la \cdot, \cdot \ra$ is a bilinear form, and
for any $x \in \E$, we have
\beq\label{Axx-tr}
\ba{rcl}
\la A x, x \ra &=& \la x x^*, A \ra.
\ea
\eeq
When $A$ is invertible, we also have
\beq\label{tr-A-Ainv}
\ba{rcl}
\la A^{-1}, A \ra &=& n.
\ea
\eeq
If the operator $W$ is positive semidefinite, and $A
\preceq A_1$ for some self-adjoint linear operator $A_1 : \E
\to \E^*$, then $\la W, A \ra \leq \la W, A_1 \ra$.
Similarly, if $A$ is positive semidefinite and $W \preceq
W_1$ for some self-adjoint linear operator $W_1 : \E^* \to
\E$, then $\la W, A \ra \leq \la W_1, A \ra$. When $A$ is
positive definite, and $R : \E \to \E^*$ is a self-adjoint
linear operator, $\la A^{-1}, R \ra$ equals the sum
of the eigenvalues of $R$ with respect to the operator $A$.
In particular, if $R$ is positive semidefinite, then all its
eigenvalues with respect to $A$ are non-negative, and the
maximal one can be bounded by the trace:
\beq\label{tr-ubd}
\ba{rcl}
R &\preceq& \la A^{-1}, R \ra \, A.
\ea
\eeq

\section{Greedy Quasi-Newton Updates}\label{sec-gr-qn}

Let $A : \E \to \E^*$ be a self-adjoint positive definite
linear operator. In this section, we consider a class of
quasi-Newton updating rules for approximating $A$.

Let $G : \E \to \E^*$ be a self-adjoint linear operator,
such that
\beq\label{qn-AG}
\ba{rcl}
A &\preceq& G,
\ea
\eeq
and let $u \in \E$ be a direction. Consider the following
family of updates, parameterized by a scalar $\tau \in \R$.
If $Gu \neq Au$, define
\beq\label{def-broyd}
\ba{rcl}
\Broyd_{\tau}(G, A, u) &\Def& \tau \left[ G -
\frac{Auu^*G + Guu^*A}{\la Au,u \ra} + \left( \frac{\la Gu,u
\ra}{\la Au,u \ra} + 1 \right) \frac{Auu^*A}{\la Au,u \ra}
\right] \\
&& + \ (1-\tau) \left[ G - \frac{(G-A)uu^*(G-A)}{\la
(G-A)u,u \ra} \right].
\ea
\eeq
Otherwise, if $Gu = Au$, define $\Broyd_{\tau}(G, A, u)
\Def G$.

Note that, for $\tau=0$, formula~\eqref{def-broyd}
corresponds to the well-known \emph{SR1 update},
\beq\label{def-sr}
\ba{rcl}
\SR(G, A, u) &\Def& G - \frac{(G-A)uu^*(G-A)}{\la (G-A)u,u
\ra},
\ea
\eeq
and, for $\tau=1$, it corresponds to the well-known
\emph{DFP update}:
\beq\label{def-dfp}
\ba{rcl}
\DFP(G, A, u) &\Def& G - \frac{Auu^*G + Guu^*A}{\la Au,u
\ra} + \left( \frac{\la Gu,u \ra}{\la Au,u \ra} + 1 \right)
\frac{Auu^*A}{\la Au,u \ra}.
\ea
\eeq
Thus, \eqref{def-broyd} describes the \emph{Broyden family}
of quasi-Newton updates (see
\cite[Section~6.3]{NocedalWright2006}), and can be written
as the linear combination of DFP and SR1
updates:\footnote{Usually, the Broyden family is defined as
the linear combination of the DFP and BFGS updates. Here we
use alternative (but equivalent) parametrization of this
family, which is more convenient for our purposes.}
$$
\ba{rcl}
\Broyd_{\tau}(G, A, u) &=& \tau \DFP(G, A, u) + (1 - \tau)
\SR(G, A, u).
\ea
$$
Our main interest will be in the class, described by the
values $\tau \in [0, 1]$, i.e. in the convex combination
of the DFP and SR1 updates. Note in particular, that this
subclass includes another well-known update---\emph{BFGS}.
Indeed, for
\beq\label{def-tau-bfgs}
\ba{rcl}
\tau_{\BFGS} &\Def& \frac{\la Au,u \ra}{\la Gu,u \ra}
\;\refIN{qn-AG}\; (0, 1),
\ea
\eeq
we have $1 - \tau_{\BFGS} = \frac{\la (G-A)u,u \ra}{\la
Gu,u \ra}$, and thus
\beq\label{def-bfgs}
\ba{rcl}
\Broyd_{\tau_{\BFGS}}(G, A, u) &=& G - \frac{\la (G-A)u,u
\ra}{\la Gu,u \ra} \frac{(G-A)uu^*(G-A)}{\la (G-A)u,u \ra}
\\ && + \ \frac{\la Au,u \ra}{\la Gu,u \ra} \left[
-\frac{Auu^*G + Guu^*A}{\la Au,u \ra} + \left( \frac{\la
Gu,u \ra}{\la Au,u \ra} + 1 \right) \frac{Auu^*A}{\la Au,u
\ra} \right] \\ &=& G - \frac{(G-A)uu^*(G-A)}{\la Gu,u \ra}
- \frac{Auu^*G+Guu^*A}{\la Gu,u \ra} \\
&& + \left( \frac{\la Gu,u \ra}{\la Au,u \ra} + 1 \right)
\frac{Auu^*A}{\la Gu,u \ra} \\ &=& G - \frac{Guu^*G}{\la
Gu,u \ra} + \frac{Auu^*A}{\la Au,u \ra} \;\Def\; \BFGS(G, A,
u).
\ea
\eeq
This is the classic BFGS formula for direction $u$.

Let us show that the Broyden family is monotonic in the
parameter $\tau$.
\BL\label{lm-mon}
If \eqref{qn-AG} holds, then, for any $u \in \E$, $\tau_1,
\tau_2 \in \R$, such that $\tau_1 \leq \tau_2$,
$$
\ba{rcl}
\Broyd_{\tau_1}(G, A, u) &\preceq& \Broyd_{\tau_2}(G, A,
u).
\ea
$$
\EL

\begin{proof}
Suppose that $Gu \neq Au$ since otherwise the claim is
trivial. Then,
$$
\ba{rcl}
&&\Broyd_{\tau}(G, A, u) \;\refEQ{def-broyd}\; G -
\frac{(G-A)uu^*(G-A)}{\la (G-A)u,u \ra} \\ &&\qquad\qquad
+ \ \tau \left[ \frac{(G-A)uu^*(G-A)}{\la (G-A)u,u \ra} -
\frac{Auu^*G+Guu^*A}{\la Au,u \ra} + \left( \frac{\la Gu,u
\ra}{\la Au,u \ra} + 1 \right) \frac{Auu^*A}{\la Au,u \ra}
\right].
\ea
$$
Denote $s \Def \frac{(G-A)u}{\la (G-A)u,u \ra} - \frac{Au}
{\la
Au,u \ra}$. Then,
$$
\ba{rl}
\la (G-A)u,u \ra s s^* &=\ \frac{(G-A)uu^*(G-A)}{\la
(G-A)u,u \ra} + \frac{\la (G-A)u,u \ra}{\la Au,u \ra}
\frac{Auu^*A}{\la Au,u \ra} - \frac{(G-A)uu^*A +
Auu^*(G-A)}{\la Au,u \ra} \\
&=\ \frac{(G-A)uu^*(G-A)}{\la
(G-A)u,u \ra} - \frac{Auu^*G + Guu^*A}{\la Au,u \ra} +
\left( \frac{\la Gu,u \ra}{\la Au,u
\ra} + 1 \right) \frac{Auu^*A}{\la Au,u \ra}.
\ea
$$
Therefore,
$$
\ba{rcl}
\Broyd_{\tau}(G, A, u) &=& G - \frac{(G-A)uu^*(G-A)}{\la
(G-A)u,u \ra} + \tau \la (G-A)u,u \ra s s^*.
\ea
$$
The claim now follows from the fact that $\la (G-A)u,u \ra s
s^* \succeq 0$ in view of \eqref{qn-AG}.
\end{proof}

Next, let us show that the relation \eqref{qn-AG} can be
preserved after applying to $G$ any update from the class of
our interest. Moreover, each update from this class does
not increase the deviation from the target operator $A$.

\BL\label{lm-pos}
If, for some $\eta \geq 1$, we have
\beq\label{pos-GA}
\ba{rclrcl}
A &\preceq& G &\preceq& \eta A,
\ea
\eeq
then, for any $u \in \E$ and any $\tau \in [0, 1]$, we
also have
\beq\label{pos-GpA}
\ba{rclrcl}
A &\preceq& \Broyd_{\tau}(G, A, u) &\preceq& \eta A.
\ea
\eeq
\EL

\begin{proof}
Denote $G_+ \Def \Broyd_{\tau}(G, A, u)$ and assume that
$Gu \neq Au$ since otherwise the claim is trivial.
Using that $\tau \geq 0$ and applying Lemma~\ref{lm-mon}, we
obtain
$$
\ba{rcl}
G_+ &\succeq& \SR(G, A, u) \;\refEQ{def-sr}\; G -
\frac{(G-A)uu^*(G-A)}{\la (G-A)u,u \ra}.
\ea
$$
Let $R \Def G-A \refSE{pos-GA} 0$, and let $I_{\E}$,
$I_{\E^*}$ be the identity operators in the spaces $\E$,
$\E^*$ respectively. Then,
$$
\ba{rcl}
G_+ - A &\succeq& R - \frac{Ruu^*R}{\la Ru,u \ra} \;=\;
\left( I_{\E^*} - \frac{Ruu^*}{\la Ru,u \ra} \right) R
\left( I_{\E} - \frac{uu^*R}{\la Ru,u \ra} \right)
\;\succeq\; 0,
\ea
$$
Thus, the first relation in \eqref{pos-GpA} is proved. To
prove the second relation, we apply Lemma~\ref{lm-mon},
using that $\tau \leq 1$, to obtain
$$
\ba{rcl}
G_+ &\preceq& \DFP(G, A, u) \;\refEQ{def-dfp}\; G + \left(
\frac{\la Gu,u \ra}{\la Au,u \ra} + 1 \right)
\frac{Auu^*A}{\la Au,u \ra} - \frac{Auu^*G + Guu^*A}{\la
Au,u \ra} \\ &=& \frac{Auu^*A}{\la Au,u \ra} + \left(
I_{\E^*} - \frac{Auu^*}{\la Au,u \ra} \right) G \left(
I_{\E} - \frac{uu^*A}{\la Au,u \ra} \right) \\
&\refPE{pos-GA}& \frac{Auu^*A}{\la Au,u \ra} + \eta \left(
I_{\E^*} - \frac{Auu^*}{\la Au,u \ra} \right) A \left(
I_{\E} - \frac{uu^*A}{\la Au,u \ra} \right) \\ &=&
\frac{Auu^*A}{\la Au,u \ra} + \eta \left( A -
\frac{Auu^*A}{\la Au,u \ra} \right) \;=\; \eta A - (\eta -
1) \frac{Auu^*A}{\la Au,u \ra} \;\preceq\; \eta A.
\ea
$$
The proof is now finished.
\end{proof}

\BR
Similar results to the one from Lemma~\ref{lm-pos} have been
known for some time in the literature for different
quasi-Newton updating formulas. For example,
in~\cite{Davidon1968} and \cite{Goldfarb1969}, it was proved
for the SR1 update that if $A \preceq G$ (respectively, $G
\preceq A$), then $A \preceq G_+$ (respectively, $G_+
\preceq A$), where $G_+$ is
the result of the corresponding update. An even stronger
property was established in \cite{Fletcher1970} for the
\emph{convex Broyden class} (composed of all convex
combinations of the BFGS and DFP updates); in particular, it
was shown that if $\eta_1 A \preceq G \preceq \eta_2 A$ for
some $0 < \eta_1
\leq 1 \leq \eta_2$, then $\eta_1 A \preceq G_+ \preceq
\eta_2 A$.
\ER

Interestingly, from Lemma~\ref{lm-mon} and
Lemma~\ref{lm-pos}, it follows that, if \eqref{qn-AG} holds,
then
$$
\ba{rclrclrclrcl}
A &\preceq& \SR(G, A, u) &\preceq& \BFGS(G, A, u) &\preceq&
\DFP(G, A, u).
\ea
$$
In other words, the approximation, produced by SR1, is
better than the one, produced by BFGS, which is in turn
better than the one, produced by DFP.

Let us now justify the efficiency of update
\eqref{def-broyd} with $\tau \in [0, 1]$ in ensuring
convergence $G \to A$. For this, we introduce the
following measure of progress:
\beq\label{def-sigma}
\ba{rcl}
\sigma_A(G) &\Def& \la A^{-1}, G - A \ra \;\refEQ
{tr-A-Ainv}\; \la A^{-1}, G \ra - n.
\ea
\eeq
Thus, $\sigma_A(G)$ is the sum of the eigenvalues of the
difference $G - A$, measured with respect to the operator
$A$. Clearly, for $G$, satisfying \eqref{qn-AG}, we have
$\sigma_A(G) \geq 0$ with $\sigma_A(G) = 0$ if and only if
$G = A$. Therefore, we need to ensure that $\sigma_A(G) \to
0$ by choosing an appropriate sequence of update directions
$u$.

First, let us estimate the decrease in the measure
$\sigma_A$ for an arbitrary direction.

\BL\label{lm-sigma-prog}
Let \eqref{qn-AG} hold. Then, for any $u \in \E$ and any
$\tau \in [0, 1]$, we have
\beq\label{sigma-prog}
\ba{rcl}
\sigma_A(G) - \sigma_A(\Broyd_{\tau}(G, A, u)) &\geq&
\frac{\la (G-A)u,u \ra}{\la Au,u \ra}.
\ea
\eeq
\EL

\begin{proof}
Denote $G_+ \Def \Broyd_{\tau}(G, A, u)$ and assume that $Gu
\neq Au$ since otherwise the claim is trivial. By
Lemma~\ref{lm-mon}, we have
$$
\ba{rcl}
G - G_+ &\succeq& G - \DFP(G, A, u) \;\refEQ{def-dfp}\;
\frac{Auu^*G + Guu^*A}{\la Au,u \ra} - \left( \frac{\la Gu,u
\ra}{\la Au,u \ra} + 1 \right) \frac{Auu^*A}{\la Au,u \ra}.
\ea
$$
Therefore,
$$
\ba{rcl}
\sigma_A(G) - \sigma_A(G_+) &\refEQ{def-sigma}& \la A^{-1},
G - G_+ \ra \;\geq\; 2 \frac{\la Gu,u \ra}{\la Au,u \ra} -
\left( \frac{\la Gu,u \ra}{\la Au,u \ra} + 1 \right) \\ &=&
\frac{\la Gu,u \ra}{\la Au,u \ra} - 1 \;=\; \frac{\la (G-A)
u, u \ra}{\la Au,u \ra}.
\ea
$$
The proof is now finished.
\end{proof}

According to Lemma~\ref{lm-sigma-prog}, the choice of the
updating direction $u$ directly influences the bound on the
decrease in the measure $\sigma_A$. Ideally, we would like
to select a direction $u$, which maximizes the right-hand
side in \eqref{sigma-prog}. However, this requires finding
an eigenvector, corresponding to the maximal eigenvalue of
$G$ with respect to $A$, which might be computationally a
difficult problem. Therefore, let us consider another
approach.

Let us fix in the space $\E$ some basis:
$$
\ba{rcl}
e_1, \dots, e_n &\in& \E.
\ea
$$
With respect to this basis, we can define the following
greedily selected direction:
\beq\label{def-gr-rule}
\ba{rcl}
\bar{u}_A(G)
&\Def& \argmax\limits_{u \in \{e_1, \dots, e_n\}} \frac{\la
(G-A)u,u \ra}{\la Au,u \ra} \;=\; \argmax\limits_{u \in
\{e_1, \dots, e_n\}} \frac{\la Gu,u \ra}{\la Au,u \ra}.
\ea
\eeq
Thus, $\bar{u}_A(G)$ is a basis vector, which maximizes the
right-hand side in \eqref{sigma-prog}. Note that for certain
choices of the basis, the computation of $\bar{u}_A(G)$
might be relatively simple. For example, if $\E = \R^n$, and
$e_1, \dots, e_n$ are coordinate directions, then the
calculation of $\bar{u}_A(G)$ requires computing only the
\emph{diagonals} of the matrix representations of the
operators $G$ and $A$. The update \eqref{def-broyd},
applying the rule \eqref{def-gr-rule}, is called the
\emph{greedy quasi-Newton update}.

Let us show that the greedy quasi-Newton update decreases
the measure $\sigma_A$ with a \emph{linear} rate. For
this, define
\beq\label{def-B}
\ba{rcl}
B &\Def& \left( \sum\limits_{i=1}^n e_i e_i^* \right)^{-1}.
\ea
\eeq
Note that $B$ is a self-adjoint positive definite linear
operator from $\E$ to $\E^*$.

\BT\label{th-gr-lin}
Let \eqref{qn-AG} hold, and let $\mu, L > 0$ be such that
\beq\label{gr-mu-L}
\ba{rclrcl}
\mu B &\preceq& A &\preceq& L B.
\ea
\eeq
Then, for any $\tau \in [0, 1]$, we have
\beq\label{gr-lin}
\ba{rcl}
\sigma_A(\Broyd_{\tau}(G, A, \bar{u}_A(G))) &\leq& \left(
1 - \frac{\mu}{n L} \right) \sigma_A(G).
\ea
\eeq
\ET

\begin{proof}
Denote $G_+ \Def \Broyd_{\tau}(G, A, \bar{u}_A(G))$, and $R
\Def G - A$. By Lemma~\ref{lm-sigma-prog},
$$
\ba{rcl}
\sigma_A(G) - \sigma_A(G_+)
&\geq& \frac{\la R\bar{u}_A(G),\bar{u}_A(G) \ra}{\la
A\bar{u}_A(G),\bar{u}_A(G) \ra} \refEQ{def-gr-rule}
\max\limits_{1 \leq i \leq n} \frac{\la R e_i, e_i \ra}{\la
A e_i, e_i \ra} \refGE{gr-mu-L} \frac{1}{L} \max\limits_
{1 \leq i \leq n} \la R e_i, e_i \ra \\
&\geq& \frac{1}{n L} \sum\limits_{i=1}^n \la R e_i, e_i \ra
\;\refEQ{Axx-tr}\; \frac{1}{n L} \sum\limits_{i=1}^n \la e_i
e_i^*, R \ra \;\refEQ{def-B}\; \frac{1}{n L} \la B^{-1}, R
\ra \\
&\refGE{gr-mu-L}& \frac{\mu}{n L} \la A^{-1}, R \ra
\;\refEQ{def-sigma}\; \frac{\mu}{n L} \sigma_A(G).
\ea
$$
The proof is now finished.
\end{proof}

\BR
A simple modification of the above proof shows that the
factor $n L$ in \eqref{gr-lin} can be improved up to
$\la B^{-1}, A \ra$. However, to simplify the future
analysis, we prefer to work directly with constant $L$.
\ER

\section{Unconstrained Quadratic Minimization}
\label{sec-quad}

Let us demonstrate how we can apply the quasi-Newton
updates, described in the previous section, for minimizing
the quadratic function
\beq\label{def-quad}
\ba{rcl}
f(x) &\Def& \frac{1}{2} \la A x, x \ra - \langle b, x
\rangle, \qquad x \in \E,
\ea
\eeq
where $A : \E \to \E^*$ is a self-adjoint positive definite
linear operator, and $b \in \E^*$.

Let $B$ be the operator, defined in \eqref{def-B}, and let
$\mu, L > 0$ be such that
\beq\label{def-mu-L}
\ba{rclrcl}
\mu B &\preceq& A &\preceq& L B.
\ea
\eeq
Thus, $\mu$ is the \emph{constant of strong convexity} of
$f$, and $L$ is the \emph{Lipschitz constant} of the
gradient of $f$, both measured with respect to the operator
$B$.

Consider the following quasi-Newton scheme:
\beq\label{quad-met-qn}
\bal \hline
\textbf{Initialization:} Choose $x_0 \in \E$. Set $G_0 = L
B$. \\ \hline \textbf{For $k \geq 0$ iterate:} \\
1. Update $x_{k+1} = x_k - G_k^{-1} \nabla f(x_k)$. \\
2. Choose $u_k \in \E$ and $\tau_k \in [0, 1]$. \\
3. Compute $G_{k+1} = \Broyd_{\tau_k}(G_k, A, u_k)$. \\
   \hline
\eal
\eeq

Note that scheme \eqref{quad-met-qn} starts with $G_0 = L
B$. Therefore, its first iteration is identical to that
one of the standard \emph{gradient method}:
$$
\ba{rcl}
x_1 = x_0 - \frac{1}{L} B^{-1} \nabla f(x_0).
\ea
$$
Also, from~\eqref{def-mu-L}, it follows that $A \preceq
G_0$. Hence, in view of Lemma~\ref{lm-pos}, we have
\beq\label{quad-GA}
\ba{rcl}
A &\preceq& G_k
\ea
\eeq
for all $k \geq 0$. In particular, all $G_k$ are positive
definite, and scheme \eqref{quad-met-qn} is well-defined.

\BR\label{quad-rm-inv}
For avoiding the $O(n^3)$ complexity for computing
$G_k^{-1} \nabla f(x_k)$, it is typical for practical
implementation of scheme \eqref{quad-met-qn} to work
directly with the inverse operators $G_k^{-1}$ (or,
alternatively, with the Cholesky decomposition of $G_k$).
Due to a low-rank structure of the updates
\eqref{def-broyd}, it is possible to compute efficiently
$G_{k+1}^{-1}$ via $G_k^{-1}$ at the cost $O(n^2)$.
\ER

To estimate the convergence rate of scheme
\eqref{quad-met-qn}, let us look at the norm of the
gradient of $f$, measured with respect to $A$:
\beq\label{quad-def-lam}
\ba{rcl}
\lambda_f(x) &\Def& \| \nabla f(x) \|_A^*
\;\refEQ{def-norms}\; \la \nabla f(x), A^{-1} \nabla f(x)
\ra^{1/2}, \qquad x \in \E.
\ea
\eeq
Note that this measure of optimality is directly related to
the functional residual. Indeed, let $x^* = A^{-1} b$ be the
minimizer of \eqref{def-quad}. Then, using Taylor's formula,
we obtain
$$
\ba{rcl}
f(x) - f^* &=& \frac{1}{2} \la A (x - x^*), x - x^* \ra
\;=\; \frac{1}{2} \la A x - b, A^{-1} (A x - b) \ra \\
&\refEQ{def-quad}& \frac{1}{2} \la \nabla f(x), A^{-1}
\nabla f(x) \ra \;\refEQ{quad-def-lam}\; \frac{1}{2}
\lambda_f^2(x).
\ea
$$

The following lemma shows how $\lambda_f$ changes after one
iteration of process \eqref{quad-met-qn}.

\BL\label{lm-quad-lambda-prog}
Let $k \geq 0$, and let $\eta_k \geq 1$ be such that
\beq\label{quad-GA-ub}
\ba{rcl}
G_k &\preceq& \eta_k A.
\ea
\eeq
Then,
$$
\ba{rcl}
\lambda_f(x_{k+1}) &\leq& \left( 1 - \frac{1}{\eta_k}
\right) \lambda_f(x_k) \;=\; \frac{\eta_k - 1}{\eta_k}
\lambda_f(x_k).
\ea
$$
\EL

\begin{proof}
Indeed,
$$
\ba{rcl}
\nabla f(x_{k+1}) &=& \nabla f(x_k) + A (x_{k+1} - x_k)
\;\refEQ{quad-met-qn}\; A (A^{-1} - G_k^{-1}) \nabla f(x_k).
\ea
$$
Therefore,
$$
\ba{rcl}
\lambda_f(x_{k+1}) &\refEQ{quad-def-lam}&
\la \nabla f(x_k), (A^{-1} - G_k^{-1}) A (A^{-1} - G_k^{-1})
\nabla f(x_k) \ra^{1/2}.
\ea
$$
Note that
$$
\ba{rclrcl}
\frac{1}{\eta_k} A^{-1} &\refPE{quad-GA-ub}& G_k^{-1}
&\refPE{quad-GA}& A^{-1}.
\ea
$$
Therefore,
\beq\label{quad-lam-prog-prel}
\ba{rclrcl}
0 &\preceq& A^{-1} - G_k^{-1} &\preceq& \left( 1 -
\frac{1}{\eta_k} \right) A^{-1}.
\ea
\eeq
Consequently,
$$
\ba{rcl}
(A^{-1} - G_k^{-1}) A (A^{-1} - G_k^{-1}) &\preceq& \left( 1
- \frac{1}{\eta_k} \right)^2 A^{-1},
\ea
$$
and
$$
\ba{rcl}
\lambda_f(x_{k+1}) &\refLE{quad-lam-prog-prel}& \left( 1 -
\frac{1}{\eta_k} \right) \la \nabla f(x_k), A^{-1} \nabla
f(x_k) \ra^{1/2} \;\refEQ{quad-def-lam}\; \left( 1 -
\frac{1}{\eta_k} \right)
\lambda_f(x_k).
\ea
$$
The proof is now finished.
\end{proof}

Thus, to estimate how fast $\lambda_f(x_k)$ converges to
zero, we need to upper bound $\eta_k$. There are two ways to
proceed, depending on the choice of directions $u_k$ in
\eqref{quad-met-qn}.

First, consider the general situation, when we do not impose
any restrictions on $u_k$. In this case, we can guarantee
that $\eta_k$ stays uniformly bounded, and $\lambda_f(x_k)
\to 0$ at a \emph{linear} rate.

\BT\label{quad-th-lin}
For all $k \geq 0$, in scheme \eqref{quad-met-qn}, we have
\beq\label{quad-G-uni-bnd}
\ba{rclrcl}
A &\preceq& G_k &\preceq& \frac{L}{\mu} A,
\ea
\eeq
and
\beq\label{quad-lam-lin}
\ba{rcl}
\lambda_f(x_k) &\leq& \left( 1 - \frac{\mu}{L} \right)^k
\lambda_f(x_0).
\ea
\eeq
\ET

\begin{proof}
Since $G_0 = L B$, in view of \eqref{def-mu-L}, we have
$$
\ba{rclrcl}
A &\preceq& G_0 &\preceq& \frac{L}{\mu} A.
\ea
$$
By Lemma~\ref{lm-pos}, this implies \eqref{quad-G-uni-bnd}.
Applying now Lemma~\ref{lm-quad-lambda-prog}, we obtain
$$
\ba{rcl}
\lambda_f(x_{k+1}) &\leq& \left( 1 - \frac{\mu}{L} \right)
\lambda_f(x_k)
\ea
$$
for all $k \geq 0$, and \eqref{quad-lam-lin} follows.
\end{proof}

Note that \eqref{quad-lam-lin} is exactly the convergence
rate of the standard gradient method. Thus, according to
Theorem~\ref{quad-th-lin}, the convergence rate of scheme
\eqref{quad-met-qn} is at least as good as that of the
gradient method.

Now assume that the directions $u_k$ in scheme
\eqref{quad-met-qn} are chosen in accordance to the greedy
strategy \eqref{def-gr-rule}. Recall that, in this case, we
can guarantee that $G_k \to A$ (Theorem~\ref{th-gr-lin}).
Therefore, we can expect faster convergence from scheme
\eqref{quad-met-qn}.

\BT\label{quad-th-super}
Suppose that, for each $k \geq 0$, we choose $u_k =
\bar{u}_A(G_k)$ in scheme \eqref{quad-met-qn}. Then, for all
$k \geq 0$, we have
\beq\label{quad-GA-lin}
\ba{rclrcl}
A &\preceq& G_k &\preceq& \left( 1 + \left( 1 - \frac{\mu}{n
L} \right)^k \frac{n L}{\mu} \right) A,
\ea
\eeq
and
\beq\label{quad-super}
\ba{rcl}
\lambda_f(x_{k+1}) &\leq& \left( 1 - \frac{\mu}{n L}
\right)^k \frac{n L}{\mu} \cdot \lambda_f(x_k).
\ea
\eeq
\ET

\begin{proof}
We already know that $A \preceq G_k$. Hence,
$$
\ba{rcl}
G_k - A &\refPE{tr-ubd}& \la A^{-1}, G_k - A \ra A
\;\refEQ{def-sigma}\; \sigma_A(G_k) A,
\ea
$$
or, equivalently,
$$
\ba{rcl}
G_k &\preceq& (1 + \sigma_A(G_k)) A.
\ea
$$
At the same time, by Theorem~\ref{th-gr-lin}, we have
$$
\ba{rcl}
\sigma_A(G_k) &\leq& \left( 1 - \frac{\mu}{n L} \right)^k
\sigma_A(G_0).
\ea
$$
Note that
$$
\ba{rcl}
\sigma_A(G_0) &\refEQ{def-sigma}& \la A^{-1}, G_0 \ra - n
\;\refLE{quad-G-uni-bnd}\; \la A^{-1}, \frac{L}{\mu} A \ra
- n \;\refEQ{tr-A-Ainv}\; n \left(\frac{L}{\mu} - 1
\right) \;\leq\; \frac{n L}{\mu}.
\ea
$$
Thus, \eqref{quad-GA-lin} is proved. Applying now
Lemma~\ref{lm-quad-lambda-prog} and using the fact that $
\frac{\eta - 1}{\eta} \leq \eta - 1$ for any $\eta \geq 1$,
we obtain \eqref{quad-super}.
\end{proof}

Theorem~\ref{quad-th-super} shows that the convergence rate
of $\lambda_f(x_k)$ is \emph{superlinear}. Let us now
combine this result with Theorem~\ref{quad-th-lin} and write
down the final efficiency estimate. Denote by $k_0 \geq 0$
the number of the first iteration, for which
\beq\label{def-k0}
\ba{rcl}
(1-\frac{\mu}{n L})^{k_0} \frac{n L}{\mu} &\leq&
\frac{1}{2}.
\ea
\eeq
Clearly, $k_0 \leq \frac{n L}{\mu} \ln \frac{2 n L}{\mu}$.
According to Theorem~\ref{th-gr-lin}, during the first $k_0$
iterations,
\beq\label{lam-k0}
\ba{rcl}
\lambda_f(x_k) &\leq& \left( 1 - \frac{\mu}{L} \right)^k
\lambda_f(x_0).
\ea
\eeq
After that, by Theorem~\ref{quad-th-super}, for all $k \geq
0$, we have
$$
\ba{rcl}
\lambda_f(x_{k_0 + k + 1}) &\refLE{quad-super}& \left( 1 -
\frac{\mu}{n L} \right)^{k_0 + k} \frac{n L}{\mu}
\lambda_f(x_{k_0 + k}) \;\refLE{def-k0}\; \left( 1 -
\frac{\mu}{n L} \right)^k
\frac{1}{2} \lambda_f(x_{k_0 + k}),
\ea
$$
or, more explicitly,
$$
\ba{rcl}
\lambda_f(x_{k_0 + k}) &\leq& \lambda_f(x_{k_0})
\prod\limits_{i=0}^{k-1} \left[ \left( 1 - \frac{\mu}{n L}
\right)^i \frac{1}{2} \right] \;=\; \left( 1 - \frac{\mu}{n
L} \right)^{\sum_{i=0}^{k-1} i} \left( \frac{1}{2} \right)^k
\lambda_f(x_{k_0}) \\ &=&
\left( 1 - \frac{\mu}{n L} \right)^{\frac{k (k-1)}{2}}
\left( \frac{1}{2} \right)^k \lambda_f(x_{k_0}) \\
&\refLE{lam-k0}& \left( 1 - \frac{\mu}{n L} \right)^{\frac{k
(k-1)}{2}} \left( \frac{1}{2} \right)^k \left( 1 -
\frac{\mu}{L} \right)^{k_0} \lambda_f(x_0).
\ea
$$
Note that the first factor in this estimate depends on the
\emph{square} of the iteration counter.

To conclude, let us mention one important property of
scheme \eqref{quad-met-qn} with greedily selected $u_k$.
It turns out that, in the particular case, when $\tau_k =
0$ for all $k \geq 0$, i.e. when scheme
\eqref{quad-met-qn} corresponds to the greedy \emph{SR1}
method, it will identify the operator $A$, and
consequently, the minimizer $x^*$ of \eqref{def-quad}, in
a \emph{finite} number of steps.

\BT
Suppose that, in scheme \eqref{quad-met-qn}, for each $k
\geq 0$, we choose $u_k = \bar{u}_A(G_k)$ and $\tau_k =
0$. Then $G_k = A$ for some $0 \leq k \leq n$.
\ET

\begin{proof}
Suppose that $R_k \Def G_k - A \neq 0$ for all $0 \leq k
\leq n$. Since $R_k \succeq 0$ (see \eqref{quad-GA}), we
must have $u_k \not\in \Ker(R_k)$ in view of
\eqref{def-gr-rule}, and
$$
\ba{rcl}
R_{k+1} &\refEQ{def-sr}& R_k - \frac{R_k u_k u_k^* R_k}{\la
R_k u_k, u_k \ra}
\ea
$$
for all $0 \leq k \leq n$. From this formula, it is easily
seen that
\BI
\II[(1)] $\Ker(R_k) \subseteq \Ker(R_{k+1})$,
\II[(2)] $u_k \in \Ker(R_{k+1})$.
\EI
Thus, the dimension of $\Ker(R_k)$ grows at least by 1 at
every iteration. In particular, the dimension of
$\Ker(R_{n+1})$ must be at least $n+1$, which is impossible,
since the operator $R_{n+1}$ acts in an $n$-dimensional
vector space.
\end{proof}

It is worth noting that for other updates, such as
\eqref{def-dfp} and \eqref{def-bfgs}, the inclusion $\Ker
(R_k) \subseteq \Ker(R_{k+1})$ is, in general, no longer
valid.

\section{Minimization of General Functions}\label{sec-gen}

Now consider a general problem of unconstrained
minimization:
\beq\label{prob-gen}
\ba{rcl}
\min\limits_{x \in \E} f(x),
\ea
\eeq
where $f : \E \to \R$ is a twice differentiable function
with positive definite Hessian. Our goal is to extend the
results, obtained in the previous section, onto the
problem~\eqref{prob-gen}, assuming that the methods can
start from a sufficiently good initial point $x_0$.

Our main assumption is that the Hessians of $f$ are close to
each other in the sense that there exists a constant $M \geq
0$, such that
\beq\label{sscf}
\ba{rcl}
\nabla^2 f(y) - \nabla^2 f(x) &\preceq& M \| y - x \|_z \,
\nabla^2 f(w)
\ea
\eeq
for all $x, y, z, w \in \E$. We call such a function $f$
\emph{strongly self-concordant}. Note that strongly
self-concordant functions form a subclass of
self-concordant functions. Indeed, let us choose a point
$x \in \E$ and a direction $h \in \E$. Then, for all $t >
0$, we have
$$
\ba{rcl}
\la [\nabla^2 f(x + t h) - \nabla^2 f(x)] h, h \ra &\leq& M
t \| h \|_x^3.
\ea
$$
Dividing this inequality by $t$ and computing the limit as
$t \downarrow 0$, we obtain
$$
\ba{rcl}
D^3 f(x)[h, h, h] &\leq& M \| h \|_x^3.
\ea
$$
for all $h \in \E$. Thus, function $f$ is self-concordant
with constant $\frac{1}{2} M$ (see
\cite{NesterovNemirovski1994,Nesterov2018Lectures}).

The main example of a strongly self-concordant function is a
strongly convex function with Lipschitz continuous Hessian.
Note however that strong self-concordancy is an
\emph{affine-invariant} property.

\BE\label{ex-sscf}
Let $C : \E \to \E^*$ be a self-adjoint positive definite
operator. Suppose there exist $\beta > 0$ and $L_2 \geq 0$,
such that the function $f$ is $\beta$-strongly convex and
its Hessian is $L_2$-Lipschitz continuous with respect to
the norm $\| \cdot \|_C$. Then $f$ is strongly
self-concordant with constant $M = \frac{L_2}{\beta^{3/2}}$.
\EE

\begin{proof}
By strong convexity of $f$, we have
\beq\label{ex-sc}
\ba{rcl}
\beta C &\preceq& \nabla^2 f(x)
\ea
\eeq
for all $x \in \E$. Therefore, using the Lipschitz
continuity of the Hessian, we obtain
$$
\ba{rcl}
\nabla^2 f(y) - \nabla^2 f(x) &\preceq& L_2 \| y - x \|_C C
\;\refEQ{def-norms}\; L_2 \la C (y - x), y - x \ra^{1/2} C
\\ &\refPE{ex-sc}& \frac{L_2}{\mu^{1/2}} \la \nabla^2 f(z)
(y - x), y - x \ra^{1/2} C \\ &\refEQ{def-norms}&
\frac{L_2}{\mu^{1/2}} \| y - x \|_z C \;\refPE{ex-sc}\;
\frac{L_2}{\mu^{3/2}} \| y - x \|_z \nabla^2 f(w)
\ea
$$
for all $x, y, z, w \in \E$.
\end{proof}

Let us establish some useful relations for strongly
self-concordant functions.

\BL\label{lm-hess-bnds}
Let $x, y \in \E$, and let $r \Def \| y - x \|_x$. Then,
\beq\label{hess-xy}
\ba{rclrcl}
\frac{\nabla^2 f(x)}{1 + M r} &\preceq& \nabla^2 f(y)
&\preceq& (1 + M r) \nabla^2 f(x).
\ea
\eeq
Also, for $J \Def \int_0^1 \nabla^2 f(x + t (y - x)) dt$
and any $v \in \{ x, y \}$, we
have
\beq\label{J-xy}
\ba{rclrcl}
\frac{\nabla^2 f(v)}{1 + \frac{M r}{2}} &\preceq& J
&\preceq& \left( 1 + \frac{M r}{2} \right) \nabla^2 f(v).
\ea
\eeq
\EL

\begin{proof}
Denote $h \Def y - x$. Taking $z = w = x$ in \eqref{sscf},
we obtain
$$
\ba{rcl}
\nabla^2 f(y) - \nabla^2 f(x) &\preceq& M r \nabla^2 f(x),
\ea
$$
which gives us the second relation in \eqref{hess-xy} after
moving $\nabla^2 f(x)$ into the right-hand side.
Interchanging now $x$ and $y$ in \eqref{sscf} and taking $z
= x$, $w = y$, we get
$$
\ba{rcl}
\nabla^2 f(x) - \nabla^2 f(y) &\preceq& M r \nabla^2 f(y),
\ea
$$
which gives us the first relation in \eqref{hess-xy} after
moving $\nabla^2 f(x)$ into the right-hand side and then
dividing by $1 + M r$.

Let us now prove \eqref{J-xy} for $v=x$ (the proof for
$v=y$ is similar). Choosing $y = x + t h$ in 
\eqref{sscf} for $t > 0$, and $w = z = x$, we obtain
$$
\ba{rcl}
\nabla^2 f(x + t h) - \nabla^2 f(x) &\preceq& M \| t h \|_x
\nabla^2 f(x) \;=\; M r t \nabla^2 f(x).
\ea
$$
This gives us the second relation in \eqref{J-xy} after
integrating for $t$ from 0 to 1 and moving $\nabla^2 f(x)$
into the right-hand side. Interchanging $x$ and $y$ in
\eqref{sscf} and taking $y = x + t h$ for $t > 0$, $z = x$,
while leaving $w$ arbitrary, we get
$$
\ba{rcl}
\nabla^2 f(x) - \nabla^2 f(x + t h) &\preceq& M \| -t h \|_x
\nabla^2 f(w) \;=\; M r t \nabla^2 f(w).
\ea
$$
Hence, by integrating for $t$ from 0 to 1, we see that
$$
\ba{rcl}
\nabla^2 f(x) - J &\preceq& \frac{M r}{2} \nabla^2 f(w).
\ea
$$
Taking now $w = x + t h$ and integrating again, we obtain
$$
\ba{rcl}
\nabla^2 f(x) - J &\preceq& \frac{M r}{2} \int_0^1 \nabla^2
f(x + t h) dt \;=\; \frac{M r}{2} J,
\ea
$$
and the first inequality in \eqref{J-xy} follows after
moving $J$ to the right-hand side and dividing by $1 +
\frac{M r}{2}$.
\end{proof}

Let us now estimate the progress of a general quasi-Newton
step. As before, for measuring the progress, we use the
\emph{local norm of the gradient}:
\beq\label{def-lam}
\ba{rcl}
\lambda_f(x) &\Def& \| \nabla f(x) \|_x^*
\;\refEQ{def-norms}\; \la \nabla f(x), \nabla^2 f(x)^{-1}
\nabla f(x) \ra^{1/2}, \qquad x \in \E.
\ea
\eeq

\BL\label{lm-lam-prog}
Let $x \in \E$, and let $G : \E \to \E^*$ be a self-adjoint
linear operator, such that
\beq\label{lam-prog-G}
\ba{rclrcl}
\nabla^2 f(x) &\preceq& G &\preceq& \eta \nabla^2 f(x)
\ea
\eeq
for some $\eta \geq 1$. Let
\beq\label{qn-step}
\ba{rcl}
x_+ &\Def& x - G^{-1} \nabla f(x),
\ea
\eeq
and let $\lambda \Def \lambda_f(x)$ be such that $M \lambda
\leq 2$. Then, $r \Def \| x_+ - x \|_x \leq \lambda$, and
$$
\ba{rcl}
\lambda_f(x_+) &\leq& \left( 1 + \frac{M \lambda}{2} \right)
\frac{\eta - 1 + \frac{M \lambda}{2}}{\eta} \lambda.
\ea
$$
\EL

\begin{proof}
Denote $J \Def \int_0^1 \nabla^2 f(x + t(x_+ - x)) dt$.
Applying Taylor's formula, we obtain
\beq\label{next-grad}
\ba{rcl}
\nabla f(x_+) &=& \nabla f(x) + J (x_+ - x)
\;\refEQ{qn-step}\; J (J^{-1} - G^{-1}) \nabla f(x).
\ea
\eeq
Note that
$$
\ba{rcl}
r &=& \| x_+ - x \|_x \;\refEQ{qn-step}\; \| G^{-1} \nabla
f(x) \|_x \;\refEQ{def-norms}\; \la \nabla f(x), G^{-1}
\nabla^2 f(x) G^{-1} \nabla f(x) \ra^{1/2} \\
&\refLE{lam-prog-G}& \la \nabla f(x), G^{-1} \nabla f(x)
\ra^{1/2} \;\refLE{lam-prog-G}\; \la \nabla f(x), \nabla^2
f(x)^{-1} \nabla f(x) \ra^{1/2} \;\refEQ{def-lam}\; \lambda.
\ea
$$
Hence, in view of Lemma~\ref{lm-hess-bnds}, we have
\beq\label{hess-via-J}
\ba{rclrclrcl}
\frac{\nabla^2 f(x)}{ 1 + \frac{M \lambda}{2} } &\preceq& J
&\preceq& \left( 1 + \frac{M \lambda}{2} \right) \nabla^2
f(x),
\qquad
J &\preceq& \left( 1 + \frac{M \lambda}{2} \right) \nabla^2
f(x_+).
\ea
\eeq
Therefore,
\beq\label{lam-prog-prel}
\ba{rcl}
\lambda_f(x_+) &\refEQ{def-lam}& \la \nabla f(x_+), \nabla^2
f(x_+)^{-1} \nabla f(x_+) \ra^{1/2} \\ &\refLE{hess-via-J}&
\sqrt{1 + \frac{M \lambda}{2}} \la
\nabla f(x_+), J^{-1} \nabla f(x_+) \ra^{1/2} \\
&\refEQ{next-grad}& \sqrt{1 + \frac{M \lambda}{2}} \la
\nabla f(x), (J^{-1} - G^{-1}) J (J^{-1} - G^{-1}) \nabla
f(x) \ra^{1/2}.
\ea
\eeq

Further,
$$
\ba{rcl}
\frac{1}{1 + \frac{M \lambda}{2}} J &\refPE{hess-via-J}&
\nabla^2 f(x) \;\refPE{lam-prog-G}\; G
\;\refPE{lam-prog-G}\; \eta
\nabla^2 f(x) \;\refPE{hess-via-J}\; \eta \left( 1 + \frac{M
\lambda}{2} \right) J.
\ea
$$
Hence,
$$
\ba{rclrcl}
\frac{1}{(1 + \frac{M \lambda}{2}) \eta} J^{-1} &\preceq&
G^{-1} &\preceq& \left( 1 + \frac{M \lambda}{2} \right)
J^{-1},
\ea
$$
and
$$
\ba{rclrcl}
-\left( 1 - \frac{1}{(1 + \frac{M \lambda}{2}) \eta} \right)
J^{-1} &\preceq& G^{-1} - J^{-1} &\preceq& \frac{M
\lambda}{2} J^{-1}.
\ea
$$
Note that
$$
\ba{rcl}
1 - \frac{1}{(1 + \frac{M \lambda}{2}) \eta} &\leq& 1 -
\frac{1 - \frac{M \lambda}{2}}{\eta} \;=\; \frac{\eta - 1 +
\frac{M \lambda}{2}}{\eta},
\ea
$$
and, since $M \lambda \leq 2$,
$$
\ba{rcl}
\frac{M \lambda}{2} &=& 1 - \left( 1 - \frac{M \lambda}{2}
\right) \;\leq\; 1 - \frac{1 - \frac{M \lambda}{2}}{\eta}
\;=\; \frac{\eta - 1 + \frac{M \lambda}{2}}{\eta}.
\ea
$$
Therefore,
$$
\ba{rclrcl}
-\frac{\eta - 1 + \frac{M \lambda}{2}}{\eta} J^{-1}
&\preceq& G^{-1} - J^{-1} &\preceq& \frac{\eta - 1 + \frac{M
\lambda}{2}}{\eta} J^{-1}.
\ea
$$
Consequently,
$$
\ba{rcl}
(G^{-1} - J^{-1}) J (G^{-1} - J^{-1}) &\preceq& \left(
\frac{\eta - 1 + \frac{M \lambda}{2}}{\eta} \right)^2
J^{-1}.
\ea
$$
and thus,
$$
\ba{rcl}
\lambda_f(x_+) &\refLE{lam-prog-prel}& \sqrt{1 + \frac{M
\lambda}{2}} \frac{\eta - 1 + \frac{M \lambda}{2}}{\eta} \la
\nabla f(x), J^{-1} \nabla f(x) \ra^{1/2} \\
&\refLE{hess-via-J}& \left( 1 + \frac{M \lambda}{2} \right)
\frac{\eta - 1 + \frac{M \lambda}{2}}{\eta} \la \nabla f(x),
\nabla^2 f(x)^{-1} \nabla f(x) \ra^{1/2} \\
&\refEQ{def-lam}& \left( 1 + \frac{M \lambda}{2} \right)
\frac{\eta - 1 + \frac{M \lambda}{2}}{\eta} \lambda.
\ea
$$
The proof is now finished.
\end{proof}

Now we need to analyze what happens with the Hessian
approximation after a quasi-Newton update. Let $G$ be the
current approximation of $\nabla^2 f(x)$, satisfying, as
usual, the condition
\beq\label{idea-hess-G}
\ba{rcl}
\nabla^2 f(x) &\preceq& G.
\ea
\eeq
Using this approximation, we can compute the new test point
$$
\ba{rcl}
x_+ &=& x - G^{-1} \nabla f(x).
\ea
$$
After that, we would like to update $G$ into a new operator
$G_+$, approximating the Hessian $\nabla^2 f(x_+)$ at the
new point and satisfying the condition
$$
\ba{rcl}
\nabla^2 f(x_+) &\preceq& G_+.
\ea
$$
A natural idea is, of course, to set
\beq\label{upd-idea}
\ba{rcl}
G_+ &=& \Broyd_{\tau}(G, \nabla^2 f(x_+), u)
\ea
\eeq
for some $u \in \E$ and $\tau \in [0, 1]$. However, we
cannot do this, since update \eqref{upd-idea} is
well-defined only when
$$
\ba{rcl}
\nabla^2 f(x_+) &\preceq& G
\ea
$$
(see Section~\ref{sec-gr-qn}), which may not be true, even
though \eqref{idea-hess-G} holds. To avoid this problem, let
us apply the following \emph{correction strategy}:
\begin{enumerate}
\II Choose some $\delta \geq 0$, and set $\tilde{G} = (1 +
\delta) G$.
\II Compute $G_+$, using \eqref{upd-idea} with $G$ replaced
by $\tilde{G}$.
\end{enumerate}
Clearly, for some value of $\delta$, the condition
$\nabla^2 f(x_+) \preceq \tilde{G}$ will be valid. If, at
the same time, this $\delta$ is sufficiently small, then
the above correction strategy should not introduce too big
error.

\BL\label{lm-op-upd}
Let $x \in \E$, and let $G : \E \to \E^*$ be a self-adjoint
linear operator, such that
\beq\label{op-upd-G}
\ba{rclrcl}
\nabla^2 f(x) &\preceq& G &\preceq& \eta \nabla^2 f(x)
\ea
\eeq
for some $\eta \geq 1$. Let $x_+ \in \E$, let $r \Def \| x_+
- x \|_x$. Then
\beq
\ba{rcl}
\tilde{G} &\Def& (1 + M r) G \;\succeq\; \nabla^2 f(x_+),
\ea
\eeq
and, for all $u \in \E$ and $\tau \in [0, 1]$, we have
$$
\ba{rclrcl}
\nabla^2 f(x_+) &\preceq& \Broyd_{\tau}(\tilde{G},
\nabla^2 f(x_+), u) &\preceq& [(1 + M r)^2 \eta] \nabla^2
f(x_+),
\ea
$$
\EL

\begin{proof}
Note that
$$
\ba{rcl}
\nabla^2 f(x_+) &\refPE{hess-xy}& (1 + M r) \nabla^2 f(x)
\;\refPE{op-upd-G}\; (1 + M r) G \;=\; \tilde{G},
\ea
$$
and,
$$
\ba{rcl}
\tilde{G} &=& (1 + M r) G
\;\refPE{op-upd-G}\; (1 + M r) \eta \nabla^2 f(x)
\;\refPE{hess-xy}\; (1 + M r)^2 \eta \nabla^2 f(x_+).
\ea
$$
Thus,
$$
\ba{rclrcl}
\nabla^2 f(x_+) &\preceq& \tilde{G} &\preceq& (1 + M r)^2
\eta \nabla^2 f(x_+),
\ea
$$
and the claim now follows from Lemma~\ref{lm-pos}.
\end{proof}

Let us now make one more assumption about the function
$f$. We assume that, with respect to the operator $B$,
defined by \eqref{def-B}, the function $f$ is
\emph{strongly convex}, and its gradient is
\emph{Lipschitz continuous}, i.e. there exist $\mu, L >
0$, such that, for all $x, y \in \E$, we have
\beq\label{hess-mu-L}
\ba{rclrcl}
\mu B &\preceq& \nabla^2 f(x) &\preceq& L B.
\ea
\eeq

\BR
In fact, for our purposes, it is enough to require that
conditions \eqref{sscf}, \eqref{hess-mu-L} hold only in a
neighborhood of a solution, but, for the sake of
simplicity, we do not do this.
\ER

We are ready to write down the scheme of our quasi-Newton
methods. For simplicity, we assume that the constants $M$
and $L$ are available.

\beq\label{met-qn}
\bal \hline
\textbf{Initialization:} Choose $x_0 \in \E$. Set $G_0 = L
B$. \\ \textbf{For $k \geq 0$ iterate:} \\
1. Update $x_{k+1} = x_k - G_k^{-1} \nabla f(x_k)$. \\
2. Compute $r_k = \| x_{k+1} - x_k \|_{x_k}$ and set
   $\tilde{G}_k = (1 + M r_k) G_k$. \\
2. Choose $u_k \in \E$ and $\tau_k \in [0, 1]$. \\
4. Compute $G_{k+1} = \Broyd_{\tau_k}(\tilde{G}_k,
   \nabla^2 f(x_{k+1}), u_k)$. \\ \hline
\eal
\eeq

\BR
Similarly to Remark~\ref{quad-rm-inv}, in a practical
implementation of scheme \eqref{met-qn}, one should work
directly with the inverse operators $G_k^{-1}$, or with
the Cholesky decomposition of $G_k$. Note that the
correction step $\tilde{G}_k = (1 + M r_k) G_k$ does not
affect the complexity of the iteration.
\ER

As before, we present two convergence results for scheme
\eqref{met-qn}. The first one establishes linear
convergence and can be seen as a generalization of
Theorem~\ref{quad-th-lin}. Note that for this result the
directions $u_k$ in the method (\ref{met-qn}) can be
chosen arbitrarily.

\BT\label{th-lin}
Suppose the initial point $x_0$ is sufficiently close
to the solution:
\beq\label{lam-ini}
\ba{rcl}
M \lambda_f(x_0) &\leq& \frac{\ln \frac{3}{2}}{4}
\frac{\mu}{L}.
\ea
\eeq
Then, for all $k \geq 0$, we have
\beq\label{op-lin}
\ba{rcl}
\nabla^2 f(x_k) &\preceq& G_k \;\preceq\; e^{2 M
\sum_{i=0}^{k-1} \lambda_f(x_i)} \frac{L}{\mu} \nabla^2
f(x_k) \;\preceq\; \frac{3 L}{2 \mu} \nabla^2 f(x_k),
\ea
\eeq
and
\beq\label{lam-lin}
\ba{rcl}
\lambda_f(x_k) &\leq& \left( 1 - \frac{\mu}{2 L} \right)^k
\lambda_f(x_0).
\ea
\eeq
\ET

\begin{proof}
In view of \eqref{hess-mu-L}, we have
$$
\ba{rclrcl}
\nabla^2 f(x_0) &\preceq& G_0 &\preceq& \frac{L}{\mu}
\nabla^2 f(x_0).
\ea
$$
Therefore, for $k=0$, both \eqref{op-lin} and
\eqref{lam-lin} are satisfied.

Now let $k \geq 0$, and suppose \eqref{op-lin},
\eqref{lam-lin} have already been proved for all $0 \leq k'
\leq k$. Denote $\lambda_k \Def
\lambda_f(x_k)$, $r_k \Def \| x_{k+1} - x_k \|_{x_k}$, and
\beq\label{def-eta}
\ba{rcl}
\eta_k &\Def& e^{2 M \sum_{i=0}^{k-1} \lambda_i}
\frac{L}{\mu}.
\ea
\eeq
Note that
\beq\label{sum-lam}
\ba{rcl}
M \sum\limits_{i=0}^k \lambda_i &\refLE{lam-lin}& M
\lambda_0 \sum\limits_{i=0}^k \left( 1 - \frac{\mu}{2 L}
\right)^i \;\leq\; \frac{2 L}{\mu} M \lambda_0
\;\refLE{lam-ini}\; \frac{\ln \frac{3}{2}}{2}.
\ea
\eeq

Applying Lemma~\ref{lm-lam-prog}, we obtain that
\beq\label{r-less-lam}
\ba{rcl}
r_k &\leq& \lambda_k
\ea
\eeq
and
\beq\label{lam-next}
\ba{rcl}
\lambda_{k+1} &\leq& \left( 1 + \frac{M \lambda_k}{2}
\right) \frac{\eta_k - 1 + \frac{M \lambda_k}{2}}{\eta_k}
\lambda_k \;=\; \left( 1 + \frac{M \lambda_k}{2} \right)
\left( 1 -
\frac{1 - \frac{M \lambda_k}{2}}{\eta_k} \right) \lambda_k.
\ea
\eeq
Note (using the inequality $1 - t \geq e^{-2t}$, valid at
least for $0 \leq t \leq \frac{1}{2}$) that
$$
\ba{rcl}
\frac{1 - \frac{M \lambda_k}{2}}{\eta_k} &\geq& e^{-M
\lambda_k} \eta_k^{-1} \;\refEQ{def-eta}\; e^{-M \lambda_k -
2 M \sum_{i=0}^{k-1}
\lambda_i} \frac{\mu}{L} \;\geq\; e^{-2 M \sum_{i=0}^k
\lambda_i} \frac{\mu}{L} \;\refGE{sum-lam}\; \frac{2 \mu}{3
L},
\ea
$$
and also (using $\ln(1+t) \leq t$, valid for $t \geq 0$)
that
$$
\ba{rcl}
\frac{M \lambda_k}{2} &\refLE{lam-ini}& \frac{\ln
\frac{3}{2}}{8} \frac{\mu}{L} \;\leq\; \frac{\mu}{16 L}.
\ea
$$
Hence,
$$
\ba{rcl}
\left( 1 + \frac{M \lambda_k}{2} \right) \left( 1 - \frac{1
- \frac{M \lambda_k}{2}}{\eta_k} \right) &\leq& \left( 1 +
\frac{\mu}{16 L} \right) \left( 1 -
\frac{2 \mu}{3 L} \right) \;\leq\; 1 - \left( \frac{2}{3} -
\frac{1}{16} \right) \frac{\mu}{L} \;\leq\; 1 - \frac{\mu}{2
L}.
\ea
$$
Consequently,
$$
\ba{rcl}
\lambda_{k+1} &\refLE{lam-next}& \left( 1 - \frac{\mu}{2 L}
\right) \lambda_k \;\refLE{lam-lin}\; \left( 1 -
\frac{\mu}{2 L} \right)^{k+1}
\lambda_0.
\ea
$$
Finally, from Lemma~\ref{lm-op-upd}, it follows that
\beq\label{G-hess-bnds}
\ba{rcl}
\nabla^2 f(x_{k+1}) &\preceq& G_{k+1}
\;\preceq\; (1 + M r_k)^2 \eta_k \nabla^2 f(x_{k+1}) \\
&\refPE{r-less-lam}& (1 + M \lambda_k)^2 \eta_k \nabla^2 f
(x_{k+1})
\;\preceq\; e^{2 M \lambda_k} \eta_k \nabla^2 f(x_{k+1}) \\
&\refEQ{def-eta}& e^{2 M \sum_{i=0}^k \lambda_i}
\frac{L}{\mu} \nabla^2 f(x_{k+1}) \;\refPE{sum-lam}\;
\frac{3 L}{2 \mu} \nabla^2 f(x_{k+1}).
\ea
\eeq
Thus, \eqref{op-lin}, \eqref{lam-lin} are valid for
$k'=k+1$, and we can continue by induction.
\end{proof}

Now let us analyze the greedy strategy. First, we analyze
how the Hessian approximation measure \eqref{def-sigma}
changes after one iteration.

\BL\label{lm-upd-sigma}
Let $x \in \E$, and let $G : \E \to \E^*$ be a self-adjoint
linear operator, such that be such that $\nabla^2 f(x)
\preceq G$. Let $x_+ \in \E$, let $r \Def \| x_+ - x \|_x$,
and let
\beq\label{def-tilde-G}
\ba{rcl}
\tilde{G} &\Def& (1 + M r) G.
\ea
\eeq
Then, for any $\tau \in [0, 1]$, we have
$$
\ba{rcl}
\sigma_{x_+}(\Broyd_{\tau}(\tilde{G}, \nabla^2 f(x_+),
\bar{u}_{x_+}(G))) &\leq& (1-\frac{\mu}{n L}) (1 + M
r)^2 \left( \sigma_x(G) + \frac{2 n M r}{1 + M r} \right).
\ea
$$
\EL

\begin{proof}
We already know from Lemma~\ref{lm-op-upd} that $\nabla^2
f(x_+) \preceq \tilde{G}$. Also note that $\bar{u}_{x_+}
(\tilde{G}) = \bar{u}_{x_+}(G)$ (see \eqref{def-gr-rule}).
Hence, by Theorem~\ref{th-gr-lin}, we have
$$
\ba{rcl}
\sigma_{x_+}(\Broyd_{\tau}(\tilde{G}, \nabla^2 f(x_+),
\bar{u}_{x_+}(G)))
\leq (1-\frac{\mu}{n L}) \sigma_{x_+}(\tilde{G}).
\ea
$$
Further,
$$
\ba{rcl}
\sigma_{x_+}(\tilde{G}) &\refEQ{def-sigma}& \la \nabla^2
f(x_+)^{-1}, \tilde{G} \ra - n \;\refEQ{def-tilde-G}\; (1
+ M r) \la \nabla^2 f(x_+)^{-1}, G \ra - n \\ &\refLE
{hess-xy}& (1 + M r)^2 \la \nabla^2 f(x)^{-1}, G \ra - n
\;\refEQ{def-sigma}\; (1 + M r)^2 \left( \sigma_x(G) + n
\right) - n \\ &=& (1 + M r)^2 \sigma_x(G) +
n ((1 + M r)^2 - 1) \\ &=& (1 + M r)^2 \sigma_x(G) + 2 n M r
\left( 1 + \frac{M r}{2} \right) \\ &\leq& (1 + M r)^2
\left( \sigma_x(G) + \frac{2 n M r}{1 + M r} \right).
\ea
$$
The proof is now finished.
\end{proof}

Now we can prove superlinear convergence. In what follows,
we assume that $n \geq 2$.

\BT\label{th-super}
Suppose that, in scheme \eqref{met-qn}, for each $k \geq 0$
we take $u_k = \bar{u}_{x_{k+1}}(G_k)$. And suppose
that the initial point $x_0$ is
sufficiently close to the solution:
\beq\label{lam-ini-super}
\ba{rcl}
M \lambda_f(x_0) &\leq& \frac{\ln 2}{4 (2 n + 1)}
\frac{\mu}{L} \quad \left( \leq \frac{\ln \frac{3}{2}}{4}
\frac{\mu}{L} \right).
\ea
\eeq
Then, for all $k \geq 0$, we have
\beq\label{op-super}
\ba{rclrcl}
\nabla^2 f(x_k) &\preceq& G_k &\preceq& \left( 1 + \left( 1
- \frac{\mu}{n L} \right)^k \frac{2 n L}{\mu} \right)
\nabla^2 f(x_k),
\ea
\eeq
and
\beq\label{lam-super}
\ba{rcl}
\lambda_f(x_{k+1}) &\leq& \left( 1 - \frac{\mu}{n L}
\right)^k \frac{2 n L}{\mu} \cdot \lambda_f(x_k).
\ea
\eeq
\ET

\begin{proof}
Denote $\lambda_k \Def \lambda_f(x_k)$ and $\sigma_k \Def
\sigma_{x_k}(G_k)$ for $k \geq 0$. In view of
Theorem~\ref{th-lin}, the first relation in \eqref{op-super}
is indeed true, and also
\beq\label{sum-lam-super}
\ba{rcl}
M \sum\limits_{i=0}^k \lambda_i &\leq& M \lambda_0
\sum\limits_{i=0}^k \left( 1 - \frac{\mu}{2 L} \right)^i
\;\leq\; \frac{2 L}{\mu} \lambda_0 \;\refLE{lam-ini-super}\;
\frac{\ln 2}{2 (2 n + 1)}.
\ea
\eeq
for all $k \geq 0$.

Let us show by induction that, for all $k \geq 0$, we have
\beq\label{sigma-p-lam}
\ba{rcl}
\sigma_k + 2 n M \lambda_k &\leq& \theta_k.
\ea
\eeq
where
\beq\label{def-theta}
\ba{rcl}
\theta_k &\Def& \left( 1 - \frac{\mu}{n L} \right)^k e^{2 (2
n + 1) M \sum_{i=0}^{k-1} \lambda_i} \frac{n L}{\mu}
\;\refLE{sum-lam-super}\; \left( 1 - \frac{\mu}{n L}
\right)^k \frac{2 n L}{\mu}.
\ea
\eeq
Indeed, since $\nabla^2 f(x_0) \preceq G_0 \preceq
\frac{L}{\mu} \nabla^2 f(x_0)$ (see \eqref{hess-mu-L}), we
have
$$
\ba{rcl}
\sigma_0 + 2 n M \lambda_0 &\refEQ{def-sigma}& \la \nabla^2
f(x_0)^{-1}, G_0 \ra - n + 2 n M \lambda_0 \\ &\leq& \la
\nabla^2 f(x_0)^{-1}, \frac{L}{\mu} \nabla^2 f (x_0) \ra -
n + 2 n M \lambda_0 \\
&\refEQ{tr-A-Ainv}& n \left( \frac{L} {\mu} - 1 \right) + 2
n M \lambda_0 \;\refLE{lam-ini-super}\; n \left( \frac{L}
{\mu} - 1 \right) + \frac{n \ln 2}{2 (2 n + 1)} \;\leq\;
\frac{n L}{\mu}.
\ea
$$
Therefore, for $k=0$, inequality \eqref{sigma-p-lam} is
satisfied. Now suppose that it is also satisfied for some $k
\geq 0$. Since $\nabla^2 f(x_k) \preceq G_k$, we know that
$$
\ba{rcl}
G_k - \nabla^2 f(x_k) &\refPE{tr-ubd}& \sigma_k \nabla^2 f
(x_k),
\ea
$$
or, equivalently,
\beq\label{G-ub-hess}
\ba{rcl}
G_k &\preceq& (1 + \sigma_k) \nabla^2 f(x_k).
\ea
\eeq
Therefore, applying Lemma~\ref{lm-lam-prog}, we obtain that
\beq\label{r-lam-super}
\ba{rcl}
r_k &\Def& \| x_{k+1} - x_k \|_{x_k} \;\leq\; \lambda_k,
\ea
\eeq
and
\beq\label{lam-next-super}
\ba{rcl}
\lambda_{k+1} &\leq& \left( 1 + \frac{M \lambda_k}{2}
\right) \frac{\sigma_k + \frac{M \lambda_k}{2}}{1 +
\sigma_k} \lambda_k \;\leq\; \left( 1 + \frac{M
\lambda_k}{2} \right) (\sigma_k
+ 2 n M \lambda_k) \lambda_k \\ &\refLE{sigma-p-lam}& \left(
1 + \frac{M \lambda_k}{2} \right) \theta_k \lambda_k
\;\leq\; e^{\frac{M \lambda_k}{2}} \theta_k \lambda_k
\;\leq\; e^{2 M \lambda_k} \theta_k \lambda_k.
\ea
\eeq
Further, by Lemma~\ref{lm-upd-sigma}, we have
$$
\ba{rcl}
\sigma_{k+1} &\leq& \left( 1 - \frac{\mu}{n L} \right) (1 +
M r_k)^2 \left( \sigma_k + \frac{2 n M r_k}{1 + M r_k}
\right) \\ &\refLE{r-lam-super}& \left( 1 - \frac{\mu}{n L}
\right) (1
+ M \lambda_k)^2 \left( \sigma_k + \frac{2 n M \lambda_k}{1
+ M \lambda_k} \right) \\ &\leq& \left( 1 - \frac{\mu}{n L}
\right) (1 + M \lambda_k)^2 (\sigma_k + 2 n M \lambda_k) \\
&\refLE{sigma-p-lam}& \left( 1 - \frac{\mu}{n L} \right) (1
+ M \lambda_k)^2 \theta_k \;\leq\; \left( 1 - \frac{\mu}{n
L} \right) e^{2 M \lambda_k} \theta_k.
\ea
$$
Note that $\frac{1}{2} \leq 1 - \frac{\mu}{n L}$ since $n
\geq 2$. Therefore,
$$
\ba{rcl}
\sigma_{k+1} + 2n M \lambda_{k+1} &\leq& \left( 1 -
\frac{\mu}{n L} \right) e^{2 M \lambda_k} \theta_k + e^{2 M
\lambda_k} \theta_k \, 2 n M \lambda_k \\ &\leq& \left( 1 -
\frac{\mu}{n L} \right) e^{2 M \lambda_k}
\theta_k + \left( 1 - \frac{\mu}{n L} \right) e^{2 M
\lambda_k} \theta_k \, 4 n M \lambda_k \\ &=& \left( 1 -
\frac{\mu}{n L} \right) e^{2 M \lambda_k} ( 1 + 4 n M
\lambda_k ) \theta_k \\ &\leq& \left( 1 - \frac{\mu}{n L}
\right) e^{2 (2n+1) M \lambda_k} \theta_k
\;\refEQ{def-theta}\; \theta_{k+1}.
\ea
$$
Thus, \eqref{sigma-p-lam} is proved.

Let us fix now some $k \geq 0$. Since $\lambda_k \geq 0$,
we have
$$
\ba{rcl}
\sigma_k &\leq& \sigma_k + 2 M \lambda_k
\;\refLE{sigma-p-lam}\; \theta_k \;\refLE{def-theta}\;
\left( 1 - \frac{\mu}{n L} \right)^k \frac{2 n L}{\mu}.
\ea
$$
This proves the second relation in \eqref{op-super} in view
of \eqref{G-ub-hess}. Finally,
$$
\ba{rl}
\lambda_{k+1} \,\refLE{lam-next-super}\, e^{2 M \lambda_k}
\theta_k \lambda_k \,\leq\, e^{2 (2 n + 1) M \lambda_k}
\theta_k \lambda_k \,\refEQ{def-theta}\,
\frac{\theta_{k+1}}{1-\frac{\mu}{n L}}
\lambda_k \,\refLE{def-theta}\, \left( 1-\frac{\mu}{n L}
\right)^k \frac{2 n L}{\mu} \lambda_k,
\ea
$$
and we obtain \eqref{lam-super}.
\end{proof}

Similarly to the quadratic case, combining
Theorem~\ref{th-lin} with Theorem~\ref{th-super}, we
obtain the following final efficiency estimate:
$$
\ba{rcl}
\lambda_f(x_{k_0 + k}) &\leq& \left( 1 - \frac{\mu}{n L}
\right)^{\frac{k (k-1)}{2}} \left( \frac{1}{2} \right)^k
\left( 1 - \frac{\mu}{2 L} \right)^{k_0} \lambda_f(x_0),
\qquad k \geq 0,
\ea
$$
where $k_0 \leq \frac{n L}{\mu} \ln \frac{2 n L}{\mu}$.

\section{Numerical Experiments}\label{sec-experim}

\subsection{Regularized Log-Sum-Exp}
\label{sec:exp-main}

In this section, we present preliminary computational
results for greedy quasi-Newton methods, applied to the
following test function\footnote{Note that we work in the
space $\E = \R^n$ and identify $\E^*$ with $\E$ in such a
way that $\la \cdot, \cdot \ra$ is the standard dot product,
and $\| \cdot \|$ is the standard Euclidean norm. Linear
operators from $\E$ to $\E^*$ are identified with $n \times
n$ matrices.}:
\beq\label{test-func}
\ba{rcl}
f(x) &\Def& \ln \left( \sum\limits_{j=1}^m e^{\la c_j, x \ra
- b_j} \right) + \frac{1}{2} \sum\limits_{j=1}^m \la c_j, x
\ra^2 + \frac{\gamma}{2} \| x \|^2, \qquad x \in \R^n,
\ea
\eeq
where $c_1, \dots, c_m \in \R^n$, $b_1, \dots, b_m \in \R$,
and $\gamma > 0$.

We compare scheme \eqref{met-qn} (which realizes GrDFP,
GrBFGS and GrSR1, depending on the choice of $\tau_k$)
with the usual gradient method (GM)\footnote{For GM, we use
the constant step size $\frac{1}{L}$, where $L$ is the
corresponding estimate of the Lipschitz constant of the
gradient, given by \eqref{test-func-L}.} and standard
quasi-Newton methods DFP, BFGS and SR1.

All the standard methods need access only to the gradient
of function $f$:
\beq\label{test-grad}
\ba{rclrcl}
\nabla f(x) &=& g(x) + \sum\limits_{j=1}^m \la c_j, x \ra
c_j + \gamma x,
\qquad
g(x) &\Def& \sum\limits_{j=1}^m \pi_j(x) c_j,
\ea
\eeq
where
$$
\ba{rcl}
\pi_j(x) &\Def& \frac{e^{\la c_j, x \ra -
b_j}}{\sum_{j'=1}^m e^{\la c_{j'}, x \ra - b_{j'}}} \;\in\;
[0, 1], \qquad j = 1, \ldots, m.
\ea
$$
Note that, for a given point $x \in \R^n$, $\nabla f(x)$ can
be computed in $O(m n)$ operations.

For greedy methods, to implement the Hessian approximation
update, at every iteration, we need to carry out some
additional operations with the Hessian
\beq\label{test-hess}
\ba{rcl}
\nabla^2 f(x) &=& \sum\limits_{j=1}^m \pi_j(x) c_j c_j^T -
g(x) g(x)^T + \sum\limits_{j=1}^m c_j c_j^T + \gamma I \\
&=& \sum\limits_{j=1}^m (\pi_j(x) + 1) c_j c_j^T - g(x)
g(x)^T + \gamma I.
\ea
\eeq
Namely, given a point $x \in \R^n$, we need to be able to
perform the following two actions:
\BI
\II For \emph{all} $1 \leq i \leq n$, compute the values
$$
\ba{rcl}
\la \nabla^2 f(x) e_i, e_i \ra &\refEQ{test-hess}&
\sum\limits_{j=1}^m (\pi_j(x) + 1) \la c_j, e_i \ra^2 - \la
g(x), e_i \ra^2 + \gamma,
\ea
$$
where $e_1, \dots, e_n$ are the basis vectors.

\II For a given direction $h \in \R^n$, compute the
Hessian-vector product
$$
\ba{rcl}
\nabla^2 f(x) h &\refEQ{test-hess}& \sum\limits_{j=1}^m
(\pi_j(x) + 1) \la c_j, h \ra c_j - \la g(x), h \ra g(x)
+ \gamma h.
\ea
$$
\EI
Let us take the basis $e_1, \ldots, e_n$, comprised of the
standard coordinate directions:
\beq\label{test-basis}
\ba{rcl}
e_i &\Def& (0, \ldots, 0, 1, 0, \ldots, 0)^T, \qquad 1 \leq
i \leq n.
\ea
\eeq
Then, both the above operations have a cost of $O(m n)$.
Thus, the cost of one iteration for all the methods under
our consideration is comparable.

Note that for basis \eqref{test-basis}, the matrix $B$,
defined by \eqref{def-B}, is the identity matrix:
$$
\ba{rcl}
B &=& I.
\ea
$$
Hence, the Lipschitz constant of the gradient of $f$ with
respect to $B$ can be taken as follows (see
\eqref{test-hess}):
\beq\label{test-func-L}
\ba{rcl}
L &=& 2 \sum\limits_{j=1}^m \| c_j \|^2 + \gamma.
\ea
\eeq
All quasi-Newton methods in our comparison start from the
same initial Hessian approximation $G_0 = L B$, and use unit
step sizes.

Finally, for greedy quasi-Newton methods, we also need to
provide an estimate of the strong self-concordancy
parameter. Note that, with respect to the operator
$\sum_{j=1}^m c_j c_j^T$, the function $f$ is 1-strongly
convex and its Hessian is 2-Lipschitz continuous (see e.g.
\cite[Ex.~1]{DoikovNesterov2019}). Hence, in view of
Example~\ref{ex-sscf}, the strong self-concordancy parameter
can be chosen as follows:
$$
\ba{rcl}
M &=& 2.
\ea
$$

The data, defining the test function \eqref{test-func}, is
randomly generated in the following way. First, we generate
a collection of random vectors
$$
\ba{rcl}
\hat{c}_1, \ldots, \hat{c}_m
\ea
$$
with entries, uniformly distributed in the interval $[-1,
1]$. Then we generate $b_1, \ldots, b_m$ from the same
distribution. Using this data, we form a preliminary
function
$$
\ba{rcl}
\hat{f}(x) &\Def& \ln \left( \sum\limits_{j=1}^m e^{\la
\hat{c}_j, x \ra - b_j} \right),
\ea
$$
and finally define
$$
\ba{rcl}
c_j &\Def& \hat{c}_j - \nabla \hat{f}(0), \qquad j = 1,
\ldots, m.
\ea
$$
Note that by construction
$$
\ba{rcl}
\nabla f(0) &\refEQ{test-grad}& \frac{1}{\sum_{j=1}^m
e^{-b_j}} \sum\limits_{j=1}^m e^{-b_j} (\hat{c}_j - \nabla
\hat{f}(0)) \;=\; 0,
\ea
$$
so the unique minimizer of our test function
\eqref{test-func} is $x^* = 0$. The starting point $x_0$ for
all methods is the same and generated randomly from the
uniform distribution on the standard Euclidean sphere of
radius $1/n$ (this choice is motivated by 
\eqref{lam-ini-super}) centered at the minimizer.

Thus, our test function \eqref{test-func} has three
parameters: the dimension $n$, the number $m$ of linear
functions, and the regularization coefficient $\gamma$. Let
us present computational results for different values of
these parameters. The termination criterion for all methods
is $f(x_k) - f(x^*) \leq \epsilon (f(x_0) - f(x^*))$.

In the tables below, for each method, we display the number
of iterations until its termination. The minus sign ($-$)
means that the method has not been able to achieve the
required accuracy after $1000n$ iterations.

\def\tmr{\multicolumn{1}{c|}{$-$}}
\def\tm{\multicolumn{1}{c}{$-$}}
\newcommand\tcclr[1]{\multicolumn{1}{|c|}{#1}}
\newcommand\tccr[1]{\multicolumn{1}{c|}{#1}}
\newcommand\tcc[1]{\multicolumn{1}{c}{#1}}
\newcommand\expn[2]{$#1 \cdot 10^{#2}$}

\begin{table}[H]\small\centering
\renewcommand\arraystretch{1}
\caption{$n = m = 50$, $\gamma = 1$}\label{tab-50-1}
\begin{tabular}{|r|r|rrr|rrr|} \hline
\tcclr{$\epsilon$} & \tccr{GM} & \tcc{DFP} & \tcc{BFGS} &
\tccr{SR1} & \tcc{GrDFP} & \tcc{GrBFGS} & \tccr{GrSR1} \\
\hline
$10^{-1}$ & 79       & 4        & 4        & 3        & 45       & 35       & 34       \\
$10^{-3}$ & 1812     & 777      & 57       & 18       & 342      & 57       & 52       \\
$10^{-5}$ & 5263     & 1866     & 107      & 29       & 738      & 72       & 58       \\
$10^{-7}$ & 8873     & 2836     & 158      & 39       & 917      & 83       & 63       \\
$10^{-9}$ & 12532    & 3911     & 203      & 48       & 1028     & 93       & 67       \\
\hline
\end{tabular}
\end{table}

\begin{table}[H]\small\centering
\renewcommand\arraystretch{1}
\caption{$n = m = 50$, $\gamma = 0.1$}\label{tab-50-01}
\begin{tabular}{|r|r|rrr|rrr|} \hline
\tcclr{$\epsilon$} & \tccr{GM} & \tcc{DFP} & \tcc{BFGS} &
\tccr{SR1} & \tcc{GrDFP} & \tcc{GrBFGS} & \tccr{GrSR1} \\
\hline
$10^{-1}$ & 76       & 4        & 4        & 3        & 44       & 33       & 33       \\
$10^{-3}$ & 2732     & 1278     & 78       & 23       & 512      & 70       & 56       \\
$10^{-5}$ & 29785    & 12923    & 254      & 57       & 3850     & 126      & 72       \\
$10^{-7}$ & \tmr     & 23245    & 346      & 74       & 6794     & 169      & 81       \\
$10^{-9}$ & \tmr     & 32441    & 381      & 79       & 8216     & 204      & 87       \\
\hline
\end{tabular}
\end{table}

\begin{table}[H]\small\centering
\renewcommand\arraystretch{1}
\caption{$n = m = 250$, $\gamma = 1$}\label{tab-250-1}
\begin{tabular}{|r|r|rrr|rrr|} \hline
\tcclr{$\epsilon$} & \tccr{GM} & \tcc{DFP} & \tcc{BFGS} &
\tccr{SR1} & \tcc{GrDFP} & \tcc{GrBFGS} & \tccr{GrSR1} \\
\hline
$10^{-1}$ & 444      & 4        & 4        & 3        & 214      & 158      & 157      \\
$10^{-3}$ & 10351    & 4743     & 98       & 21       & 3321     & 264      & 251      \\
$10^{-5}$ & 73685    & 31468    & 288      & 55       & 15637    & 350      & 274      \\
$10^{-7}$ & 159391   & 58138    & 450      & 82       & 21953    & 413      & 296      \\
$10^{-9}$ & 249492   & 85218    & 627      & 110      & 25500    & 464      & 314      \\
\hline
\end{tabular}
\end{table}

\begin{table}[H]\small\centering
\renewcommand\arraystretch{1}
\caption{$n = m = 250$, $\gamma = 0.1$}\label{tab-250-01}
\begin{tabular}{|r|r|rrr|rrr|} \hline
\tcclr{$\epsilon$} & \tccr{GM} & \tcc{DFP} & \tcc{BFGS} &
\tccr{SR1} & \tcc{GrDFP} & \tcc{GrBFGS} & \tccr{GrSR1} \\
\hline
$10^{-1}$ & 442      & 4        & 4        & 3        & 209      & 155      & 155      \\
$10^{-3}$ & 9312     & 4175     & 91       & 21       & 2686     & 258      & 251      \\
$10^{-5}$ & 207978   & 102972   & 488      & 87       & 60461    & 556      & 346      \\
$10^{-7}$ & \tmr     & \tm      & 1003     & 170      & 147076   & 792      & 391      \\
$10^{-9}$ & \tmr     & \tm      & 1407     & 233      & 212100   & 976      & 419      \\
\hline
\end{tabular}
\end{table}

We see that all quasi-Newton methods outperform the gradient
method and demonstrate superlinear convergence (from some
moment, the difference in the number of iterations between
successive rows in the table becomes smaller and smaller).
Among quasi-Newton methods (both the standard and the greedy
ones), SR1 is always better than BFGS, while DFP is
significantly worst than the other two. At the first few
iterations, the greedy methods loose to the standard ones,
but later they catch up. However, the classical SR1 method
always remains the best. Nevertheless, the greedy methods
are quite competitive.

Now let us look at the quality of Hessian approximations,
produced by the quasi-Newton methods. In the tables below,
we display the desired accuracy $\epsilon$ vs the final
Hessian approximation error (defined as the operator norm of
$G_k - \nabla^2 f(x_k)$, measured with respect to $\nabla^2
f(x_k)$). We look at the same problems as in
Table~\ref{tab-50-1} and Table~\ref{tab-250-1}.

\begin{table}[H]\small\centering
\renewcommand\arraystretch{1}
\caption{$n = m = 50$, $\gamma = 1$}
\begin{tabular}{|r|rrr|rrr|} \hline
\tcclr{$\epsilon$} & \tcc{DFP} & \tcc{BFGS} & \tccr{SR1} &
\tcc{GrDFP} & \tcc{GrBFGS} & \tccr{GrSR1} \\ \hline
$10^{-0}$ & \expn{1.6}{3}    & \expn{1.6}{3}    & \expn{1.6}{3}    & \expn{1.6}{3}    & \expn{1.6}{3}    & \expn{1.6}{3} \\
$10^{-1}$ & \expn{1.6}{3}    & \expn{1.6}{3}    & \expn{1.6}{3}    & \expn{2.7}{3}    & \expn{1.5}{3}    & \expn{1.5}{3} \\
$10^{-3}$ & \expn{1.6}{3}    & \expn{1.6}{3}    & \expn{1.6}{3}    & \expn{1.2}{3}    & \expn{1.2}{1}    & \expn{3.8}{0} \\
$10^{-5}$ & \expn{1.6}{3}    & \expn{1.6}{3}    & \expn{1.6}{3}    & \expn{2.1}{2}    & \expn{7.2}{0}    & \expn{2.6}{0} \\
$10^{-7}$ & \expn{1.6}{3}    & \expn{1.6}{3}    & \expn{1.6}{3}    & \expn{9.1}{1}    & \expn{5.6}{0}    & \expn{2.2}{0} \\
$10^{-9}$ & \expn{1.6}{3}    & \expn{1.6}{3}    & \expn{1.6}{3}    & \expn{5.2}{1}    & \expn{4.1}{0}    & \expn{1.8}{0} \\
\hline
\end{tabular}
\end{table}

\begin{table}[H]\small\centering
\renewcommand\arraystretch{1}
\caption{$n = m = 250$, $\gamma = 1$}
\begin{tabular}{|r|rrr|rrr|} \hline
\tcclr{$\epsilon$} & \tcc{DFP} & \tcc{BFGS} & \tccr{SR1} &
\tcc{GrDFP} & \tcc{GrBFGS} & \tccr{GrSR1} \\ \hline
$10^{-0}$  & \expn{4.1}{4}    & \expn{4.1}{4}    & \expn{4.1}{4}    & \expn{4.1}{4}    & \expn{4.1}{4}  & \expn{4.1}{4} \\
$10^{-1}$  & \expn{4.1}{4}    & \expn{4.1}{4}    & \expn{4.1}{4}    & \expn{7.1}{4}    & \expn{3.8}{4}  & \expn{3.9}{4} \\
$10^{-3}$  & \expn{4.1}{4}    & \expn{4.1}{4}    & \expn{4.1}{4}    & \expn{6.8}{4}    & \expn{6.6}{1}  & \expn{1.7}{1} \\
$10^{-5}$  & \expn{4.1}{4}    & \expn{4.1}{4}    & \expn{4.1}{4}    & \expn{9.4}{3}    & \expn{3.7}{1}  & \expn{1.2}{1} \\
$10^{-7}$  & \expn{4.1}{4}    & \expn{4.1}{4}    & \expn{4.1}{4}    & \expn{3.1}{3}    & \expn{2.8}{1}  & \expn{9.7}{0}    \\
$10^{-9}$  & \expn{4.1}{4}    & \expn{4.1}{4}    & \expn{4.1}{4}    & \expn{1.7}{3}    & \expn{2.2}{1}  & \expn{7.3}{0}    \\
\hline
\end{tabular}
\end{table}

As we can see from these tables, for standard quasi-Newton
methods the Hessian approximation error always stays at the
initial level. In contrast, for the greedy ones, it
decreases relatively fast (especially for GrBFGS and GrSR1).
Note also that sometimes the initial residual slightly
increases at the first several iterations (which is
noticeable only for GrDFP). This happens due to the fact
that the objective function is non-quadratic, and we apply
the correction strategy.

Note that in all the above tests we have used the same
values for the parameters $n$ and $m$. Let us briefly
illustrate what happens when, for example, $m > n$.

\begin{table}[H]\small\centering
\renewcommand\arraystretch{1}
\caption{$n = 50$, $m = 100$, $\gamma = 0.1$}
\begin{tabular}{|r|r|rrr|rrr|} \hline
\tcclr{$\epsilon$} & \tccr{GM} & \tcc{DFP} & \tcc{BFGS} &
\tccr{SR1} & \tcc{GrDFP} & \tcc{GrBFGS} & \tccr{GrSR1} \\
\hline
$10^{-1}$    & 84       & 4        & 4        & 3        &
46 & 37       & 37       \\
$10^{-3}$    & 897      & 316      & 32       & 11       &
183      & 53       & 52       \\
$10^{-5}$    & 2421     & 833      & 67       & 19       &
334      & 63       & 58       \\
$10^{-7}$    & 4087     & 1304     & 98       & 25       &
423      & 71       & 62       \\
$10^{-9}$    & 5810     & 1859     & 132      & 32       &
473      & 78       & 66       \\
\hline
\end{tabular}
\end{table}

\begin{table}[H]\small\centering
\renewcommand\arraystretch{1}
\caption{$n = 50$, $m = 200$, $\gamma = 0.1$}
\begin{tabular}{|r|r|rrr|rrr|} \hline
\tcclr{$\epsilon$} & \tccr{GM} & \tcc{DFP} & \tcc{BFGS} &
\tccr{SR1} & \tcc{GrDFP} & \tcc{GrBFGS} & \tccr{GrSR1} \\
\hline
$10^{-1}$    & 108      & 4        & 4        & 3        &
45 & 46       & 46       \\
$10^{-3}$    & 479      & 101      & 17       & 7        &
97       & 53       & 52       \\
$10^{-5}$    & 1059     & 338      & 39       & 12       &
154      & 62       & 59       \\
$10^{-7}$    & 1817     & 615      & 62       & 18       &
206      & 67       & 64       \\
$10^{-9}$    & 2659     & 807      & 81       & 21       &
234      & 73       & 68       \\
\hline
\end{tabular}
\end{table}

Comparing these tables with Table~\ref{tab-50-01}, we see
that, with the increase of $m$, all the methods generally
terminate faster. However, the overall picture is still the
same as before. The results for $m < n$ are similar, so we
do not include them.

Finally, let us present the results for the
\emph{randomized} version of scheme \eqref{met-qn}, in
which, at every step, we select the update direction
uniformly at random from the standard Euclidean sphere:
\beq\label{u-rand}
\ba{rcl}
u_k &\sim& \Unif(\mathcal{S}^{n-1}),
\ea
\eeq
where $\mathcal{S}^{n-1} \Def \{ x \in \R^n : \| x \| = 1
\}$. We call the corresponding methods RaDFP, RaBFGS and
RaSR1.

\begin{center}
\begin{minipage}{.45\textwidth}
\begin{table}[H]\small\centering
\renewcommand\arraystretch{1}
\caption{$n = m = 50$, $\gamma = 1$}
\begin{tabular}{|r|rrr|} \hline
\tcclr{$\epsilon$} & \tcc{RaDFP} & \tcc{RaBFGS} &
\tccr{RaSR1} \\
\hline
$10^{-1}$    & 35       & 29       & 34       \\
$10^{-3}$    & 566      & 102      & 64       \\
$10^{-5}$    & 1156     & 125      & 77       \\
$10^{-7}$    & 1481     & 142      & 85       \\
$10^{-9}$    & 1698     & 156      & 91       \\
\hline
\end{tabular}
\end{table}
\end{minipage}
\begin{minipage}{.45\textwidth}
\begin{table}[H]\small\centering
\renewcommand\arraystretch{1}
\caption{$n = m = 250$, $\gamma = 1$}
\begin{tabular}{|r|rrr|} \hline
\tcclr{$\epsilon$} & \tcc{RaDFP} & \tcc{RaBFGS} &
\tccr{RaSR1} \\
\hline
$10^{-1}$    & 261      & 144      & 158      \\
$10^{-3}$    & 4276     & 366      & 287      \\
$10^{-5}$    & 19594    & 517      & 346      \\
$10^{-7}$    & 33293    & 619      & 376      \\
$10^{-9}$    & 41177    & 698      & 396      \\
\hline
\end{tabular}
\end{table}
\end{minipage}
\end{center}
\vspace{1em}
It is instructive to compare these tables with
Table~\ref{tab-50-1} and Table~\ref{tab-250-1}, which
contain the results for the greedy methods on the same
problems. We see that the randomized methods are slightly
slower than the greedy ones. However, the difference is not
really significant, and, what is especially interesting, the
randomized methods do not loose superlinear convergence.

\subsection{Logistic Regression}

Now let us consider another test function, namely
\emph{$l_2$-regularized logistic regression}, which is
popular in the field of machine learning:
\beq\label{def-log-reg}
\ba{rcl}
f(x) &\Def& \sum\limits_{j=1}^m \ln(1 + e^{-b_j \la c_j, x
\ra}) + \frac{\gamma}{2} \| x \|^2, \qquad x \in \R^n,
\ea
\eeq
where $c_1, \dots, c_m \in \R^n$, $b_1, \ldots, b_m \in \{
-1, 1\}$, and $\gamma > 0$.

Note that the structure of the function \eqref{def-log-reg}
is similar to the one of \eqref{test-func}. In particular,
both the diagonal of the Hessian and the Hessian-vector
product for this function can be computed with the similar
complexity of that for computing the gradient. Also it can
be shown that the Lipschitz constant of the gradient of $f$
can be chosen in accordance with \eqref{test-func-L} but
with the coefficient $\frac{1}{4}$ instead of $2$.

We follow the same experiment design as before with only a
couple of differences. First, instead of generating the
data, defining the function \eqref{def-log-reg},
artificially, now we take it from the LIBSVM collection of
real-world data sets for binary classification
problems\footnote{The original labels $b_i$ in the
\emph{mushrooms} data set are ``$1$'' and ``$2$'' instead of
``$1$'' and ``$-1$''. Therefore, we renamed in advance the
class label ``$2$'' into ``$-1$''.} \cite{ChangLin2011}.
Second, we have found it better in practice not to apply the
correction strategy in the greedy methods (i.e. simply set
$\tilde{G}_k = G_k$ in scheme \eqref{met-qn}). This is the
only heuristic that we use. For the regularization
coefficient, we always use the value $\gamma=1$, which is a
standard choice.

Let us look at the results.

\begin{table}[H]\small\centering
\renewcommand\arraystretch{1}
\caption{Data set \emph{ijcnn1} ($n=22$, $m=49990$)}
\label{tab-ijcnn1}
\begin{tabular}{|r|r|rrr|rrr|} \hline
\tcclr{$\epsilon$} & \tccr{GM} & \tcc{DFP} & \tcc{BFGS} &
\tccr{SR1} & \tcc{GrDFP} & \tcc{GrBFGS} & \tccr{GrSR1} \\
\hline
$10^{-1}$    & 246      & 43       & 8        & 6        &
25 & 19       & 18       \\
$10^{-3}$    & 1925     & 672      & 45       & 16       &
71       & 25       & 23       \\
$10^{-5}$    & 5123     & 2007     & 85       & 25       &
145      & 32       & 23       \\
$10^{-7}$    & 8966     & 2738     & 102      & 29       &
192      & 38       & 23       \\
$10^{-9}$    & 12815    & 3269     & 118      & 33       &
215      & 43       & 24       \\
\hline
\end{tabular}
\end{table}

\begin{table}[H]\small\centering
\renewcommand\arraystretch{1}
\caption{Data set \emph{mushrooms} ($n=112$, $m=8124$)}
\label{tab-mushrooms}
\begin{tabular}{|r|r|rrr|rrr|} \hline
\tcclr{$\epsilon$} & \tccr{GM} & \tcc{DFP} & \tcc{BFGS} &
\tccr{SR1} & \tcc{GrDFP} & \tcc{GrBFGS} & \tccr{GrSR1} \\
\hline
$10^{-1}$    & 4644     & 936      & 15       & 6        &
230 & 83       & 82       \\
$10^{-3}$    & 77103    & 30594    & 105      & 24       &
1185     & 149      & 113      \\
$10^{-5}$    & \tmr       & 58221    & 166      & 34       &
1700     & 170      & 113      \\
$10^{-7}$    & \tmr       & 83740    & 217      & 42       &
1945     & 182      & 113      \\
$10^{-9}$    & \tmr       & 107471   & 257      & 48       &
2088     & 194      & 114      \\
\hline
\end{tabular}
\end{table}

\begin{table}[H]\small\centering
\renewcommand\arraystretch{1}
\caption{Data set \emph{a9a} ($n=123$, $m=32561$)}
\label{tab-a9a}
\begin{tabular}{|r|r|rrr|rrr|} \hline
\tcclr{$\epsilon$} & \tccr{GM} & \tcc{DFP} & \tcc{BFGS} &
\tccr{SR1} & \tcc{GrDFP} & \tcc{GrBFGS} & \tccr{GrSR1} \\
\hline
$10^{-1}$    & 160      & 32       & 10       & 6        &
110 & 81       & 81       \\
$10^{-3}$    & 18690    & 9229     & 145      & 38       &
2203     & 127      & 117      \\
$10^{-5}$    & \tmr       & 79014    & 411      & 88       &
23715    & 316      & 123      \\
$10^{-7}$    & \tmr       & \tm       & 553      & 113     
&
35700    & 441      & 124      \\
$10^{-9}$    & \tmr       & \tm       & 581      & 118     
&
38285    & 475      & 124      \\
\hline
\end{tabular}
\end{table}

\begin{table}[H]\small\centering
\renewcommand\arraystretch{1}
\caption{Data set \emph{w8a} ($n=300$, $m=49749$)}
\label{tab-w8a}
\begin{tabular}{|r|r|rrr|rrr|} \hline
\tcclr{$\epsilon$} & \tccr{GM} & \tcc{DFP} & \tcc{BFGS} &
\tccr{SR1} & \tcc{GrDFP} & \tcc{GrBFGS} & \tccr{GrSR1} \\
\hline
$10^{-1}$    & 10148    & 3531     & 35       & 10       &
694 & 300      & 300      \\
$10^{-3}$    & 194813   & 86315    & 178      & 34       &
1426     & 307      & 301      \\
$10^{-5}$    & \tmr       & 188561   & 300      & 54      
& 1849     & 327      & 301      \\
$10^{-7}$    & \tmr       & 255224   & 387      & 68       &
2036     & 339      & 301      \\
$10^{-9}$    & \tmr       & 264346   & 399      & 69       &
2057     & 340      & 301      \\
\hline
\end{tabular}
\end{table}

As we can see, the general picture is the same as for the
previous test function. In particular, the DFP update is
always much worse than BFGS and SR1. The greedy methods are
competitive with the standard ones and often outperform them
for high values of accuracy.

\section{Discussion}

We have presented the greedy quasi-Newton methods, that are
based on the updating formulas from the Broyden family and
use greedily selected basis vectors for updating Hessian
approximations. For these methods, we have established
explicit non-asymptotic rate of local superlinear
convergence for the iterates and also a linear convergence
for the deviations of Hessian approximations from the
correct Hessians.

Clearly, there is a number of open questions. First, at
every iteration, our methods need to compute the greedily
selected basis vector. This requires additional information
beyond just the gradient of the objective function (such as
the diagonal of the Hessian). However, many problems, that
arise in applications, possess certain structure (separable,
sparse, etc.), for which the corresponding computations have
a cost similar to that of the gradient evaluation (such as
the test function in our experiments). Nevertheless, ideally
it is desirable to get rid of the necessity in this
auxiliary information at all. A natural idea might be to
replace the greedy strategy with a \emph{randomized} one.
Indeed, as can be seen from our experiments, the
corresponding scheme \eqref{met-qn}, \eqref{u-rand}
demonstrate almost the same performance as the greedy one.
Therefore, one can expect that it should be possible to
establish similar theoretical results about its
superlinear convergence. Nevertheless, at the moment, we do
not know how to do this. Although it is not difficult to
show that, in terms of \emph{expectations}, the randomized
strategy still preserves the linear convergence of Hessian
approximations (see \cite{GowerRichtarik2017}), it is not
clear how to proceed after this in proving the superlinear
convergence of the iterates, even in the quadratic case. The
main difficulty, arising in the analysis, is that, at
some moment, one needs to take the expectation of the
product of random variables with known expectations, but the
random variables themselves are \emph{non-independent}.

Second, we have analyzed together a whole class of Hessian
approximation updates by essentially upper bounding all its
members via the worst one---DFP. Thus, all the efficiency
guarantees, that we have established, might be too
pessimistic for other members of this class such as BFGS,
and especially SR1. Indeed, in our experiments, we have seen
that the convergence properties of these three methods might
differ quite significantly. It is therefore desirable to
refine our current analysis and obtain separate estimates
for different updates.

Third, note that our current results do not prove anything
about the rate of superlinear convergence of the standard
quasi-Newton methods. Of course, it would be interesting to
obtain the corresponding estimates and compare them to the
ones, that we have established in this work.

Finally, apart from the quadratic case, we have not
addressed at all the question of global convergence.

In any case, we believe that the ideas and the theoretical
analysis, presented in this paper, will be useful for future
advances in the theory of quasi-Newton methods.

\section*{Acknowledgments}
The authors would like to thank two anonymous referees
for their useful comments and suggestions.


\begin{thebibliography}{10}

\bibitem{Davidon1959}
W. Davidon. Variable metric method for minimization. Argonne
   National Laboratory Research and Development Report 5990
   (1959).

\bibitem{FletcherPowell1963}
R. Fletcher and M. Powell. A rapidly convergent descent
   method for minimization. \emph{Computer Journal},
   \textbf{6}(2), 163-168 (1963).

\bibitem{Broyden1967}
C. Broyden. Quasi-Newton methods and their application to
   function minimization. \emph{Mathematics of Computation},
   \textbf{21}(99), 368-381 (1967).

\bibitem{Davidon1968}
W. Davidon. Variance algorithm for minimization.
   \emph{Computer Journal}, \textbf{10}(4), 406-410 (1968).

\bibitem{Goldfarb1969}
D. Goldfarb. Sufficient conditions for the convergence of a
   variable metric algorithm. \emph{Optimization}, ed. R.
   Fletcher, 273-281, Academic Press, London (1969).

\bibitem{Broyden1970p1}
C. Broyden. The convergence of a class of double-rank
   minimization algorithms: 1.~General considerations.
   \emph{IMA Journal of Applied Mathematics}, \textbf{6}(1),
   76-90 (1970).

\bibitem{Broyden1970p2}
C. Broyden. The convergence of a class of double-rank
   minimization algorithms: 2.~The new algorithm. \emph{IMA
   Journal of Applied Mathematics}, \textbf{6}(3), 222-231
   (1970).

\bibitem{Fletcher1970}
R. Fletcher. A new approach to variable metric algorithms.
   \emph{Computer Journal}, \textbf{13}(3), 317-322 (1970).

\bibitem{Goldfarb1970}
D. Goldfarb. A family of variable-metric methods derived by
   variational means. \emph{Mathematics of Computation},
   \textbf{24}(109), 23-26 (1970).

\bibitem{Shanno1970}
D. Shanno. Conditioning of quasi-Newton methods for function
   minimization. \emph{Mathematics of Computation},
   \textbf{24}(111), 647-656 (1970).

\bibitem{Powell1971}
M. Powell. On the convergence of the variable metric
   algorithm. \emph{IMA Journal of Applied Mathematics},
   \textbf{7}(1), 21-36 (1971).

\bibitem{BroydenDennisMore1973}
C. Broyden, J. Dennis, and J. Mor{\'e}. On the local and
   superlinear convergence of quasi-Newton methods.
   \emph{IMA Journal of Applied Mathematics},
   \textbf{12}(3), 223-245 (1973).

\bibitem{DennisMore1974}
J. Dennis and J. Mor{\'e}. A characterization of superlinear
   convergence and its application to quasi-Newton methods.
   \emph{Mathematics of Computation}, \textbf{28}(126),
   549-560 (1974).

\bibitem{DennisMore1977}
J. Dennis and J. Mor{\'e}. Quasi-Newton methods, motivation
   and theory. \emph{SIAM Review}, \textbf{19}(1), 46-89
   (1977).

\bibitem{PshenichnyiDanilin1978}
B. Pshenichny{\u{i}} and I. Danilin. Numerical methods in
extremal problems. \emph{Mir Publishers} (1978).

\bibitem{Stachurski1981}
A. Stachurski. Superlinear convergence of Broyden's bounded
   $\theta$-class of methods. \emph{Mathematical
   Programming}, \textbf{20}(1), 196-212 (1981).

\bibitem{GriewankToint1982}
A. Griewank and P. Toint. Local convergence analysis for
   partitioned quasi-Newton updates. \emph{Numerische
   Mathematik}, \textbf{39}(3), 429-448 (1982).

\bibitem{ByrdNocedalYuan1987}
R. Byrd, J. Nocedal, and Y. Yuan. Global convergence of a
   class of quasi-Newton methods on convex problems.
   \emph{SIAM Journal on Numerical Analysis},
   \textbf{24}(5), 1171-1190 (1987).

\bibitem{ByrdNocedal1989}
R. Byrd and J. Nocedal. A tool for the analysis of
   quasi-Newton methods with application to unconstrained
   minimization. \emph{SIAM Journal on Numerical Analysis},
   \textbf{26}(3), 727-739 (1989).

\bibitem{ConnGouldToint1991}
A. Conn, N. Gould, and P. Toint. Convergence of quasi-Newton
   matrices generated by the symmetric rank one update.
   \emph{Mathematical Programming}, \textbf{50}(1-3),
   177-195 (1991).

\bibitem{EngelsMartinez1991}
J. Engels and H. Mart{\'i}nez. Local and superlinear
   convergence for partially known quasi-Newton methods.
   SIAM Journal on Optimization, \textbf{1}(1), 42-56
   (1991).

\bibitem{NesterovNemirovski1994}
Y. Nesterov and A. Nemirovskii. Interior-point polynomial
   algorithms in convex programming. \emph{SIAM},
   \textbf{13} (1994).

\bibitem{YabeYamaki1996}
H. Yabe and N. Yamaki. Local and superlinear convergence of
   structured quasi-Newton methods for nonlinear
   optimization. \emph{Journal of the Operations Research
   Society of Japan}, \textbf{39}(4), 541-557 (1996).

\bibitem{WeiYuYuanLian2004}
Z. Wei, G. Yu, G. Yuan, and Z. Lian. The superlinear
   convergence of a modified BFGS-type method for
   unconstrained optimization. \emph{Computational
   Optimization and Applications}, \textbf{29}(3), 315-332
   (2004).

\bibitem{NocedalWright2006}
J. Nocedal and S. Wright. Numerical optimization.
   \emph{Springer Science \& Business Media} (2006).

\bibitem{YabeOgasawaraYoshino2007}
H. Yabe, H. Ogasawara, and M. Yoshino. Local and superlinear
   convergence of quasi-Newton methods based on modified
   secant conditions. \emph{Journal of Computational and
   Applied Mathematics}, \textbf{205}(1), 617-632 (2007).

\bibitem{ChangLin2011}
C. Chang and C. Lin. LIBSVM: A library for support vector
   machines. \emph{ACM transactions on intelligent systems
   and technology (TIST)}, \textbf{2}(3), 1-27 (2011).

\bibitem{LewisOverton2013}
A. Lewis and M. Overton. Nonsmooth optimization via
   quasi-Newton methods. \emph{Mathematical Programming},
   \textbf{141}(1-2), 135-163 (2013).

\bibitem{GowerGoldfarbRichtarik2016}
R. Gower, D. Goldfarb, and P. Richt{\'a}rik. Stochastic
   block BFGS: squeezing more curvature out of data.
   \emph{International Conference on Machine Learning},
   1869-1878 (2016).

\bibitem{GowerRichtarik2017}
R. Gower and P. Richt{\'a}rik. Randomized quasi-Newton
   updates are linearly convergent matrix inversion
   algorithms. \emph{SIAM Journal on Matrix Analysis and
   Applications}, \textbf{38}(4), 1380-1409 (2017).

\bibitem{MokhtariEisenRibeiro2018}
A. Mokhtari, M. Eisen, and A. Ribeiro. IQN: An incremental
   quasi-Newton method with local superlinear convergence
   rate. \emph{SIAM Journal on Optimization},
   \textbf{28}(2), 1670-1698 (2018).

\bibitem{Nesterov2018Lectures}
Y. Nesterov. Lectures on convex optimization.
   \emph{Springer}, \textbf{137} (2018).

\bibitem{DoikovNesterov2019}
N. Doikov and Y. Nesterov. Minimizing uniformly convex
   functions by cubic regularization of Newton method.
   \emph{arXiv}, 1905.02671 (2019)

\bibitem{SunTranDinh2019}
T. Sun and Q. Tran-Dinh. Generalized self-concordant
   functions: a recipe for Newton-type methods.
   \emph{Mathematical Programming}, \textbf{178}, 145-213
   (2019).

\bibitem{GaoGoldfarb2019}
W. Gao and D. Goldfarb. Quasi-Newton methods: superlinear
   convergence without line searches for self-concordant
   functions. \emph{Optimization Methods and Software},
   \textbf{34}(1), 194-217 (2019).

\bibitem{KovalevEtAl2020}
D. Kovalev, R. Gower, P. Richt{\'a}rik, and A. Rogozin. Fast
   linear convergence of randomized BFGS. \emph{arXiv
   preprint} arXiv:2002.11337 (2020).
\end{thebibliography}
\end{document}